\documentclass[a4paper,10pt,reqno]{amsart}
\usepackage{amssymb,latexsym,amsmath,amsfonts}
\usepackage{amsthm}
\usepackage[mathscr]{eucal}
\usepackage{rotating}
\usepackage{multirow}
\numberwithin{equation}{section}
\newtheorem{theorem}{Theorem}

\newtheorem{proposition}{Proposition}

\theoremstyle{definition}

\newtheorem{example}[theorem]{Example}

\theoremstyle{remark}
\newtheorem{remark}{Remark}

\newcommand {\sgn}{\mbox{sgn}}

\begin{document}

\title[Decay Rate Classification of Regularly Varying ODEs]
{Classification of convergence rates of solutions of perturbed
ordinary differential equations with regularly varying
nonlinearity}

\author{John A. D. Appleby}
\address{School of Mathematical Sciences, Dublin City University, Glasnevin, Dublin 9, Ireland}
\email{john.appleby@dcu.ie} \urladdr{webpages.dcu.ie/\textasciitilde
applebyj}

\author{Denis D. Patterson}
\address{School of Mathematical
Sciences, Dublin City University, Glasnevin, Dublin 9, Ireland}
\email{denis.patterson2@mail.dcu.ie}

\thanks{John Appleby gratefully acknowledges Science Foundation
Ireland for the support of this research under the Mathematics
Initiative 2007 grant 07/MI/008 ``Edgeworth Centre for Financial
Mathematics''. Denis Patterson is supported by the Irish Research Council grant GOIPG/2013/402 ``Persistent and strong dependence in growth rates of solutions of stochastic and deterministic
functional differential equations with applications to finance''.} 
\date{24 December 2015}

\begin{abstract}
In this paper we consider the rate of convergence of solutions of a scalar ordinary differential equation which is a perturbed 
version of an autonomous equation with a globally stable equilibrium. Under weak assumptions on the nonlinear mean reverting force, 
we demonstrate that the convergence rate is preserved when the perturbation decays more rapidly than a critical rate. At the critical rate, 
the convergence to equilibrium is slightly slower than the unperturbed equation, and when the perturbation decays more slowly than the critical rate, 
the convergence to equilibrium is strictly slower than that seen in the unperturbed equation. In the last case, under strengthened assumptions, a new convergence rate is recorded which depends on the convergence rate of the perturbation. The latter result relies on the function being regularly varying at the equilibrium with index greater than unity; therefore, for this class of regularly varying problems, a classification of the convergence rates is obtained.  
\end{abstract}

\maketitle

\section{Introduction}
In this paper we classify the rates of convergence to a limit of the solutions of a scalar ordinary differential equation
\begin{equation} \label{eq.odepert}
x'(t)=-f(x(t))+g(t), \quad t>0; \quad x(0)=\xi.
\end{equation}
We assume that the unperturbed equation
\begin{equation} \label{eq.ode}
y'(t)=-f(y(t)), \quad t>0; \quad y(0)=\zeta
\end{equation}
has a unique globally stable equilibrium (which we set to be at zero). This is characterised by the condition
\begin{equation} \label{eq.fglobalstable}
xf(x)>0 \quad\text{for $x\neq 0$,} \quad f(0)=0.
\end{equation}
In order to ensure that both \eqref{eq.ode} and \eqref{eq.odepert} have
continuous solutions, we assume
\begin{equation} \label{eq.fgcns}
  f\in C(\mathbb{R};\mathbb{R}), \quad g\in C([0,\infty);\mathbb{R}).
\end{equation}
The condition \eqref{eq.fglobalstable} ensures that any solution of \eqref{eq.odepert} is \emph{global} i.e., that
\[
\tau:=\inf\{t>0\,:\,x(t)\not\in (-\infty,\infty)\}=+\infty.
\]
We also ensure that there is exactly one continuous solution of both \eqref{eq.odepert} and \eqref{eq.ode} by assuming
\begin{equation} \label{eq.floclip}
\text{$f$ is locally Lipschitz continuous on $\mathbb{R}$}.
\end{equation}
In \eqref{eq.ode} or \eqref{eq.odepert}, we assume that $f(x)$ does not have linear leading order behaviour as $x\to 0$; 
moreover, we do not ask that $f$ forces solutions of \eqref{eq.ode} to hit zero in finite time. Since $f$ is continuous, we are free to define
\begin{equation} \label{def.F}
F(x)=\int_x^1 \frac{1}{f(u)}\,du, \quad x>0,
\end{equation}
and avoiding solutions of \eqref{eq.ode} to hitting zero in finite time forces
\begin{equation} \label{eq.Ftoinfty}
\lim_{x\to 0^+} F(x)=+\infty.
\end{equation}
We notice that $F:(0,\infty)\to\mathbb{R}$ is a strictly decreasing function, so it has an inverse $F^{-1}$. Clearly, \eqref{eq.Ftoinfty} implies that
\[
\lim_{t\to\infty} F^{-1}(t)=0.
\]
The significance of the functions $F$ and $F^{-1}$ is that they enable us to determine the rate of convergence of solutions of \eqref{eq.ode} to zero, because
$F(y(t))-F(\zeta)=t$ for $t\geq 0$ or $y(t)=F^{-1}(t+F(\zeta))$ for $t\geq 0$. It is then of interest to ask whether solutions of \eqref{eq.odepert} will still converge to zero as $t\to\infty$, and how this convergence rate modifies according to the asymptotic behaviour of $g$.

In order to do this with reasonable generality we find it convenient and natural to assume at various points that the functions $f$ and $g$ are regularly varying.
We recall that a measurable function $f:(0,\infty)\to (0,\infty)$ with $f(x)>0$ for $x>0$ is said to be regularly varying at $0$ with index $\beta\in\mathbb{R}$ if 
\[
\lim_{x\to 0^+} \frac{f(\lambda x)}{f(x)}=\lambda^\beta, \quad \text{for all $\lambda>0$}.
\]
A measurable function $h:[0,\infty)\to [0,\infty)$ with $h(t)>0$ for $t\geq 0$ is said to regularly varying at infinity with index $\alpha\in\mathbb{R}$ if 
\[
\lim_{t\to\infty} \frac{h(\lambda t)}{h(t)}=\lambda^\alpha, \quad \text{for all $\lambda>0$}.
\]
We use the notation $f\in \text{RV}_0(\beta)$ and $h\in\text{RV}_\infty(\alpha)$. Many useful properties of regularly varying functions, including those employed here, are recorded in Bingham, Goldie and Teugels~\cite{BGT}.  

The main result of the paper, which \emph{characterises} the rate of convergence of solutions of \eqref{eq.odepert} to zero, can be summarised as follows: suppose that $f$ is regularly varying at zero with index $\beta>1$, and that $g$ is positive and
regularly varying at infinity, in such a manner that
\[
\lim_{t\to\infty} \frac{g(t)}{(f\circ F^{-1})(t)}=:L\in [0,\infty]
\]
exists. If $L=0$, the solution of \eqref{eq.odepert} inherits the rate of decay to zero of $y$, in the sense that
\[
\lim_{t\to\infty} \frac{F(x(t))}{t}=1.
\]
If $L\in (0,\infty)$ we can show that the rate of decay to zero is slightly slower, obeying
\[
\lim_{t\to\infty} \frac{F(x(t))}{t}=\Lambda=\Lambda(L)\in (0,1)
\]
and a formula for $\Lambda$ purely in terms of $L$ and $\beta$ can be found. Finally, in the case that $L=+\infty$ it can be shown that
\[
\lim_{t\to\infty} \frac{F(x(t))}{t}=0,
\]
If it is presumed that $g$ is regularly varying at infinity with negative index, or $g$ is slowly varying and is asymptotic to a decreasing function, then the exact rate of convergence can be found, namely that $\lim_{t\to\infty} f(x(t))/g(t)=1$. These asymptotic results are proven by constructing appropriate upper and lower solutions to the differential equation \eqref{eq.odepert} as in Appleby and Buckwar~\cite{appbuck}.

In some cases, we do not need the full strength of the regular variation hypotheses: when $L=0$, all that is needed is the asymptotic monotonicity of $f$ close to zero; on the other hand, the hypothesis $\beta>1$ seems to be important in the case when $L\in (0,\infty]$. 
If $f$ is regularly varying with index $\beta=1$, examples exist for which $L=+\infty$, but $F(x(t))/t\to 1$ as $t\to\infty$. Therefore the conditions under which 
this asymptotic characterisation holds seem best suited to the case when $f$ is regularly varying at $0$ with index $\beta>1$. 

There is a wealth of literature concerning the use of regular variation in analysing the asymptotic behaviour of ordinary differential equations, and the field is very active. Besides work of Avakumovi\'c in 1947 on equations of Thomas--Fermi type in~\cite{Akumovic}, some of the earliest work is due to Mari\'c and Tomi\'c \cite{MarTom:76, MarTom:80} concerning the asymptotic behaviour of nonlinear second order ordinary differential equations, with linear second order equations 
being treated in depth by Omey~\cite{Omey:81}. An important monograph summarising themes in the research up to the year 2000 is Maric~\cite{Maric2000}. More recently highly nonlinear and nonautonomous second--order differential equations of Emden--Fowler type have been studied with regularly varying state--dependence and non--autonomous multiplier, ~\cite{KusMan:2011a,KusMan:2011b, Matucci2014,Matucci2015}, as well as solutions of nonautonomous linear functional differential equations with time--varying delay \cite{KusanoMaric:2006a} and higher--order differential equations \cite{KusanoMaric:2006b}. 
Another important strand of research on the exact asymptotic behaviour of non--autonomous ordinary differential equations (of first and higher order) in which the equations have regularly varying coefficients has been developed. For recent contributions, see for example work of Evtukhov and co--workers (e.g., Evtukhov and Samoilenko~\cite{EvSam:2011}) and Koz\'ma~\cite{Kozma:2012}, as well as the references in these papers. These papers tend to be concerned with non--autonomous features which are \emph{multipliers} of the regularly--varying state dependent terms, in contrast to the presence of the nonautonomous term $g$ in \eqref{eq.odepert}, which might be thought of as \emph{additive}. Despite this extensive literature and active research concerning regular variation and asymptotic behaviour of ordinary differential equations, and despite the fact that our analysis deals with first--order equations only, it would appear that the results presented in this work are new.

One of the motivations for this work is to consider the asymptotic behaviour of solutions of stochastic differential equations of It\^o type with state--independent diffusion coefficient in which the drift function is $-f$ and $f$ is regularly varying. In Appleby and Patterson~\cite{appleby_patterson} we have developed some of the results in the present paper to allow solutions to change sign and  impose integral rather than pointwise conditions on the forcing term to preserve decay rates to equilibrium. Such extensions are crucial in providing a comprehensive treatment of SDEs of the type mentioned above. A further motivation for the current work is to extend results in~\cite{appleby_patterson} to deal with SDEs with slowly decaying diffusion coefficient, and the results presented here which deal with slowly decaying $g$ should form an important ingredient in performing this analysis. 

The paper is organised as follows: in Section 2 the main results of the paper are discussed, and notation and supporting results outlined. Section 3 contains examples showing the scope of the theorem. Some of these examples show why the conditions of the main results are difficult to relax without fundamentally altering the asymptotic behaviour of solutions. The proofs of the main results are given in the final Section 4.

\section{Mathematical Preliminaries, Discussion of Hypotheses and Statement of Main Results}
In this section we introduce some common notation and list known properties of regular, slow and rapidly varying functions. We also discuss the hypotheses used in 
the paper, and then lay out and discuss the main results of the paper. 
\subsection{Notation and properties of regularly varying functions}
Throughout the paper, the set of real numbers is denoted by $\mathbb{R}$. We let $C(I;J)$ stand for the space of continuous functions which map $I$ onto $J$, where $I$ and $J$ are typically intervals in $\mathbb{R}$. Similarly, the space of differentiable functions with continuous derivative mapping $I$ onto $J$ is denoted by $C^1(I;J)$. If $h$ and $j$ are real--valued functions defined on $(0,\infty)$ and $\lim_{t\to\infty} h(t)/j(t)=1$, we sometimes use the standard asymptotic notation $h(t)\sim j(t)$ as $t\to\infty$. Similarly, if $h$ and $j$ obey $\lim_{t\to 0^+} h(t)/j(t)=1$, we write $h(t)\sim j(t)$ as $t\to\infty$.

The results quoted in this short section concerning regularly varying functions at infinity may all be found in Chapter 1 in~\cite{BGT}. They are listed below for the completeness of the exposition. Properties listed below of functions that are regularly varying at 0 may be deduced from properties of functions which are regularly varying at infinity by exploiting the fact that if $f\in \text{RV}_0(\beta)$, then $h:(0,\infty)\to (0,\infty)$ defined by 
\[
h(t)=\frac{1}{f(1/t)}, \quad t>0
\]
is in $\text{RV}_\infty(\beta)$.
\begin{itemize}
\item[(i)] Composition and reciprocals:
If $h\in \text{RV}_\infty(-\theta)$ for $\theta\geq 0$ and $h(t)\to 0$ as $t\to\infty$, and $\phi\in \text{RV}_0(\beta)$ for $\beta>0$, then 
$\phi\circ h\in \text{RV}_\infty(-\theta\beta)$. If $h\in \text{RV}_\infty(\theta)$, then $1/h\in \text{RV}_\infty(-\theta)$, while $\phi\in \text{RV}_0(\beta)$
implies $1/\phi\in \text{RV}_0(-\beta)$.
\item[(ii)] Inverses:  
If there is $\eta<0$ such that $\phi\in \text{RV}_0(\eta)$ (so that $\phi(x)\to \infty$ as $x\to 0^+$) and $\phi:(0,\delta)\to (0,\infty)$ is invertible, then $\phi^{-1}\in \text{RV}_\infty(1/\eta)$. If there is $\eta>0$ such that $\phi\in \text{RV}_0(\eta)$ (so that $\phi(x)\to 0$ as $x\to 0^+$) and $\phi:(0,\delta)\to (0,\infty)$ is invertible, then $\phi^{-1}\in \text{RV}_0(1/\eta)$. Similarly, if there is $\theta>0$ such that $h\in \text{RV}_\infty(-\theta)$ (so that $h(t)\to 0$ as $t\to\infty$) and $h:(T,\infty)\to (0,\infty)$ is invertible, then $h^{-1} \in \text{RV}_0(-1/\theta)$.
\item[(iii)] Preservation of asymptotic order:  
If $x,y\in C([0,\infty);(0,\infty))$ are such that $\lim_{t\to\infty}x(t)=\lim_{t\to\infty} y(t)=0$, and $x(t)/y(t)\to 1$ as $t\to\infty$, and $\phi\in \text{RV}_0(\beta)$ for $\beta\neq 0$, then 
\[
\lim_{t\to\infty} \frac{\phi(x(t))}{\phi(y(t))}=1.
\]
Similarly if $x,y\in C([0,\infty);(0,\infty))$ are such that $\lim_{t\to\infty}x(t)=+\infty$, $\lim_{t\to\infty} y(t)=+\infty$, and $x(t)/y(t)\to 1$ as $t\to\infty$, and 
$h\in \text{RV}_\infty(\theta)$ for $\theta\neq 0$, then 
\[
\lim_{t\to\infty} \frac{h(x(t))}{h(y(t))}=1.
\]
\item[(iv)] Integration:  
If $\phi$ in $\text{RV}_0(\beta)$ for $\beta>1$, then 
\[
\lim_{x\to 0^+} \frac{\int_{x}^1 1/\phi(u)\,du}{\frac{1}{\beta-1}\frac{x}{\phi(x)}}=1.
\]
\item[(v)] Smooth approximation:  
If $h$ is in $\text{RV}_\infty(-\theta)$ for $\theta>0$, then there exists $j\in C^1((0,\infty);(0,\infty))$ which is also in $\text{RV}_\infty(-\theta)$ such that $j'(t)<0$ for all $t>0$ and 
\[
\lim_{t\to\infty} \frac{h(t)}{j(t)}=1, \quad \lim_{t\to\infty}\frac{tj'(t)}{j(t)}=-\theta.
\] 
Similarly, if $\phi\in \text{RV}_0(\beta)$ for $\beta>0$, then there exists $\varphi\in C^1((0,\infty),\infty))\cap \text{RV}_0(\beta)$ such that $\varphi'(x)>0$ for all $x>0$ and
\[
\lim_{x\to 0^+} \frac{\phi(x)}{\varphi(x)}=1, \quad \lim_{x\to 0^+} \frac{x\varphi'(x)}{\varphi(x)}=\beta.
\]
A slightly weaker result holds for slowly varying functions at $\infty$: if $h$ is in $\text{RV}_\infty(0)$, 
then there exists $j\in C^1((0,\infty);(0,\infty))$ which is also in $\text{RV}_\infty(0)$ such that  
\[
\lim_{t\to\infty} \frac{h(t)}{j(t)}=1, \quad \lim_{t\to\infty}\frac{tj'(t)}{j(t)}=0.
\] 
It is part of e.g., Theorem 1.3.3 in~\cite{BGT}.
\item[(vi)] Uniform asymptotic behaviour on compact intervals: 
We observe that if $h\in \text{RV}_\infty(-\theta)$, then for any $c>0$ we have 
\[
\lim_{t\to\infty} \frac{h(t-c)}{h(t)}=1.
\]
\end{itemize}
Some further terminology should be introduced. We say that a function $\phi$ is \textit{slowly varying at $0$} if $\phi\in \text{RV}_0(0)$ and that a function $h$ is \textit{slowly varying at infinity} if $h\in\text{RV}_\infty(0)$. A function $h:(0,\infty)\to (0,\infty)$ is said to be \textit{rapidly varying of index $\infty$ at infinity} if 
\[
\lim_{t\to\infty} \frac{h(\lambda t)}{h(t)}=
\left\{
\begin{array}{cc}
0, & \lambda<1, \\
1, & \lambda=1, \\
+\infty, & \lambda>1.
\end{array}
\right.
\] 
For such a function $h$ we write $h\in \text{RV}_\infty(\infty)$. Analogously, a function $h:(0,\infty)\to (0,\infty)$ is said to be \textit{rapidly varying of index $-\infty$ at infinity} if 
\[
\lim_{t\to\infty} \frac{h(\lambda t)}{h(t)}=
\left\{
\begin{array}{cc}
+\infty, & \lambda<1, \\
1, & \lambda=1, \\
0, & \lambda>1.
\end{array}
\right.
\] 
For such a function $h$ we write $h\in \text{RV}_\infty(-\infty)$. Together, these two classes of functions are described as being rapidly varying at infinity.
We can extend naturally this notation to deal with rapid variation at zero.  
Suppose that $\phi:(0,\infty)\to (0,\infty)$ is measurable such that 
\[
\lim_{x\to 0^+} \frac{\phi(\lambda x)}{\phi(x)}
=
\left\{
\begin{array}{cc}
+\infty, & \lambda>1, \\
1, & \lambda=1, \\
0, & \lambda<1.
\end{array}
\right.
\]
We write $\phi\in \text{RV}_0(\infty)$. On the other hand, if 
\[
\lim_{x\to 0^+} \frac{\phi(\lambda x)}{\phi(x)}
=
\left\{
\begin{array}{cc}
0, & \lambda>1, \\
1, & \lambda=1, \\
+\infty, & \lambda<1,
\end{array}
\right.
\]
we write $\phi\in \text{RV}_0(-\infty)$. There is a connection between rapidly and slowly varying functions through inverses. It is a fact that if $h\in \text{RV}_\infty(\infty)$ (which forces $h(t)\to \infty$ as $t\to\infty$) and $h$ is invertible, then $h^{-1}\in \text{RV}_\infty(0)$. 
\subsection{Discussion of hypotheses}
In order to simplify the analysis in this paper, we assume that
\begin{equation} \label{eq.gposxipos}
g(t)>0 \quad t>0; \quad x(0)=\xi>0.
\end{equation}
This has the effect of restricting the solutions of \eqref{eq.odepert} to be positive for all $t\geq 0$ and assists us in characterising convergence rates according to the rate of decay of $g$. We will show in further work that this sign assumption can be lifted, and our desired asymptotic characterisation will be for the most part preserved. Moreover, it transpires that the results in this work can be used to prove results when the sign restriction is relaxed, by means of comparison proofs.
Our asymptotic results also tacitly assume that $g(t)\to 0$ as $t\to\infty$, but in further work we show that this assumption can also be relaxed, while maintaining results on the rate of decay of solutions of \eqref{eq.odepert}. In fact, as mentioned above the analysis in this paper will enable the almost sure rate of convergence rates of solutions of (It\^o) stochastic differential equations with state independent noise intensity to be analysed. 

The results of this paper can rapidly be extended in the case that $g(t)<0$ for all $t\geq 0$ and $\xi<0$. In this case, consider $x_-(t)=-x(t)$ for $t\geq 0$, 
$g_-(t)=-g(t)$ for $t\geq 0$, $\xi_-=-\xi$ and $f_-(x)=-f(-x)$ for $x\in \mathbb{R}$. Then 
\[
x_-'(t)=f(x(t))-g(t)=-f_-(x_-(t))+g_-(t), \quad t>0; \quad x_-(0)=\xi_-.
\]
Clearly, $g_-$ and $\xi_-$ now obey \eqref{eq.gposxipos} and $f_-$ still obeys \eqref{eq.fglobalstable}, \eqref{eq.floclip}, and if $g$ is continuous so is $g_-$.
Therefore, we can prove asymptotic results for $x_-$ using the results given in this paper, and therefore readily recover those results for $x$.

Any discussion of convergence rates of $x(t)\to 0$ as $t\to\infty$ implicitly assumes that the desired convergence actually occurs. Rather than making additional
assumptions on $f$ and $g$ in this paper which guarantee convergence, we will assume that the convergence occurs. One result which guarantees that $x(t)\to 0$ as $t\to\infty$ is nonetheless recorded below, because additional hypotheses on $g$ follow from our assumptions in many cases.
\begin{proposition} \label{prop.xto0ginL1}
Suppose that $f$ obeys \eqref{eq.fglobalstable}, that $f$ and $g$ obey \eqref{eq.fgcns}, and that $g\in L^1(0,\infty)$.
Then every continuous solution $x$ of \eqref{eq.odepert} obeys
\begin{equation} \label{eq.xto0}
\lim_{t\to\infty} x(t)=0.
\end{equation}
\end{proposition}
In the case when $g$ is not integrable, but $g(t)\to 0$ as $t\to\infty$, it can be shown that either \eqref{eq.xto0} holds or $x(t)\to\infty$
as $t\to\infty$ (see e.g.,~\cite{JAJC:2011szeged}). Solutions of \eqref{eq.xto0} exhibit a type of local stability: if the initial condition $\xi$ and sup norm of $g$ are sufficiently small, \eqref{eq.xto0} is true. A sufficient condition which rules out unbounded solutions, and therefore guarantees \eqref{eq.xto0} for all initial conditions, is
\begin{equation} \label{eq.fphiinfty}
\liminf_{x\to\infty} f(x)>0.
\end{equation}
See~\cite{JAJGAR:2009} for example. In the case when $f(x)\to 0$ as $x\to\infty$, the relationship between the rate of decay of $g(t)\to 0$ as $t\to\infty$ and $f(x)\to 0$ as $x\to\infty$ becomes important: for a given $f$, if the rate of decay of $g$ is too slow and the initial condition is too large, then $x(t)\to\infty$ as $t\to\infty$. However, if $g$ decays more quickly than a certain rate, it can be shown that \eqref{eq.xto0} holds for all initial conditions. Moreover, under some additional hypotheses,
a critical rate of decay of $g$ can be identified, in the sense that if $g$ decays more slowly to zero than this rate, solutions can escape to infinity, while 
if it decays faster than the critical rate, solutions obey \eqref{eq.xto0} for all initial conditions. For further details, we refer the reader to~\cite{JAJC:2011szeged} and the references therein. It is interesting to note that a condition of the form 
\eqref{eq.fphiinfty} is unnecessary for almost sure asymptotic stability in SDEs, and accordingly, this hypothesis is not appealed to in 
~\cite{appleby_patterson}.

\subsection{Main results}
We now state and discuss our results precisely. In our first result, we can show that the global convergence of solutions of \eqref{eq.odepert}, as well as the rate of convergence of solutions to $0$ is preserved provided the perturbation $g$ decays sufficiently rapidly. In order to guarantee this, we request only that 
$f$ be asymptotic to a monotone function close to zero: no regular  variation is needed.
\begin{theorem} \label{theorem.gsmall}
Suppose that $f$ obeys \eqref{eq.fglobalstable}, \eqref{eq.floclip} and that $F$ defined by \eqref{def.F} obeys \eqref{eq.Ftoinfty}.
Suppose further that $f$ and $g$ obey \eqref{eq.fgcns} and that \eqref{eq.gposxipos} holds. Suppose that there exists $\phi$ such that
\begin{equation} \label{eq.fasyphimon}
\lim_{x\to 0^+} \frac{f(x)}{\phi(x)}=1, \quad \text{ $\phi$ is increasing on $(0,\delta)$}.
\end{equation}
If
\begin{equation} \label{eq.gfF0}
\lim_{t\to\infty} \frac{g(t)}{(f\circ F^{-1})(t)}=0,
\end{equation}
then the unique continuous solution of \eqref{eq.odepert} obeys
\begin{equation}  \label{eq.Fxtto1}
\lim_{t\to\infty} \frac{F(x(t))}{t}=1.
\end{equation}
\end{theorem}
Immediately Theorem~\ref{theorem.gsmall} presents a question: is it possible to find slower rates of decay of $g(t)\to 0$ as $t\to\infty$ than exhibited in \eqref{eq.gfF0}, for which the solution $x$ of \eqref{eq.odepert} still decays at the rate of the unperturbed equation, as characterised by \eqref{eq.Fxtto1}?
In some sense, our next theorem says that the rate of decay of $g$ in \eqref{eq.gfF0} \emph{cannot} be relaxed, at least for functions $f$ which are regularly varying at $0$ with index $\beta>1$, or which are rapidly varying at zero. 

In the case when $f$ is regularly varying at $0$ with index 1 (and $f(x)/x\to 0$ 
as $x\to 0$), the condition \eqref{eq.gfF0} is not necessary in order to preserve the rate of decay embodied by \eqref{eq.Fxtto1}. This claim is confirmed by 
the following example. It also suggests, in the case when $f$ is regularly varying at zero with index $1$, that a more careful analysis is needed to characterise
the asymptotic behaviour of solutions of \eqref{eq.odepert}. 
\begin{example}
Suppose $\delta\in (e^{-(\sqrt{2}-1)},1)$ and define $f(x)=x/\log(1/x)$ for $x\in (0,\delta)$ and let $f(0)=0$. We see that $f\in \text{RV}_0(1)$ 
and $f(x)/x\to 0$ as $x\to 0^+$. Suppose that 
\[
g(t)=e^{-\sqrt{2}(1+t)^{1/2}+(1+t)^{1/3}}\frac{1}{(1+t)^{2/3}}\cdot \frac{\frac{5\sqrt{2}}{6}-\frac{1}{3}(1+t)^{-1/6}}{\sqrt{2}-(1+t)^{-1/6}}, \quad t\geq 0.
\]
Then $g$ is continuous and $g(t)>0$ for all $t\geq 0$. Consider the initial value problem 
\[
x'(t)=-f(x(t))+g(t), \quad t>0; \quad x(0)=e^{-(\sqrt{2}-1)}. 
\]
It can be verified $x(t)=\exp(-\sqrt{2}(1+t)^{1/2}+(1+t)^{1/3})$ for $t\geq 0$ satisfies this initial value problem, and is therefore its unique continuous solution. Defining $F(x)=\int_x^\delta 1/f(u)\,du$ for $x\in (0,\delta)$ we see that 
\begin{gather*}
F(x)=\frac{1}{2}\left(\log^2(1/x)-\log^2(1/\delta)\right), \quad x\in (0,\delta), \\ 
F^{-1}(x)=\exp\left(-\sqrt{2x+\log^2(1/\delta)}\right), \quad x>0.
\end{gather*}
Hence $f\circ F^{-1}$ is well--defined on $[0,\infty)$, and we can rapidly show that 
\[
\lim_{t\to\infty} \frac{(f\circ F^{-1})(t)}{e^{-\sqrt{2}t^{1/2}}\frac{1}{t^{1/2}\sqrt{2}}}=1.
\]
Therefore, it follows that 
\[
\lim_{t\to\infty} \frac{g(t)}{(f\circ F^{-1})(t)}=+\infty.
\]
Since formulae for $F$ and $x$ are known, it is easily checked that $F(x(t))/t\to 1$ as $t\to\infty$. Therefore, it can be seen that \eqref{eq.gfF0}
is violated, $f$ is regularly varying (with index $1$) at $0$, and all other hypotheses of Theorem~\ref{theorem.gsmall} are satisfied, but nonetheless the solution of the initial value problem \eqref{eq.odepert} obeys \eqref{eq.Fxtto1}.
\end{example}

We now turn to asking how the rate of decay changes when \eqref{eq.gfF0} is relaxed, and $f$ is regularly varying at $0$ with index $\beta>1$ or is rapidly varying at zero. 
\begin{theorem} \label{theorem.FxtLambda}
Suppose that $f$ obeys \eqref{eq.fglobalstable}, \eqref{eq.floclip} and that $F$ defined by \eqref{def.F} obeys \eqref{eq.Ftoinfty}.
Suppose further that $f$ and $g$ obey \eqref{eq.fgcns} and that \eqref{eq.gposxipos} holds. 
Let $x$ be the unique continuous solution of \eqref{eq.odepert}. 
Suppose that there exists $\phi$ such that \eqref{eq.fasyphimon} holds, and suppose further that there exists $L>0$ such that
\begin{equation} \label{eq.gfFL}
\lim_{t\to\infty} \frac{g(t)}{(f\circ F^{-1})(t)}=L.
\end{equation}
Then $x(t)\to 0$ as $t\to\infty$. 
\begin{itemize}
\item[(i)] If $f\in \text{RV}_0(\beta)$ for $\beta>1$, then  
\begin{equation}  \label{eq.limFxtLambda} 
\lim_{t\to\infty} \frac{F(x(t))}{t}=\Lambda_\ast(L)\in (0,1),
\end{equation}
where $\Lambda_\ast$ is the unique solution of $(1-\Lambda_\ast)\Lambda_\ast^{-\beta/(\beta-1)}=L$.
\item[(ii)] If $f\circ F^{-1}\in \text{RV}_\infty(-1)$ and $F^{-1}\in \text{RV}_\infty(0)$, then 
\[
\lim_{t\to\infty} \frac{F(x(t))}{t}=\Lambda_\ast(L)\in (0,1),
\]
where $\Lambda_\ast$ is the unique solution of $(1-\Lambda_\ast)\Lambda_\ast^{-1}=L$, or $\Lambda_\ast=1/(L+1)$.
\end{itemize}
\end{theorem}
If $y$ is the solution of \eqref{eq.ode}, we have that $y(t)/F^{-1}(t)\to 1$ as $t\to\infty$. Moreover, in the case when $\beta>1$, as $F^{-1}\in \text{RV}_\infty(-1/(\beta-1))$,
we have 
\[
\lim_{t\to\infty} \frac{x(t)}{y(t)}
=\lim_{t\to\infty}\frac{x(t)}{F^{-1}(t)}=\lim_{t\to\infty}\frac{F^{-1}(\Lambda_\ast t)}{F^{-1}(t)} = \Lambda_\ast^{-1/(\beta-1)}>1.
\]
Therefore, the solution of \eqref{eq.odepert} is of the same order as the solution of \eqref{eq.ode}, but decays more slowly by a factor depending on $L$.
In the second case, when $F^{-1}\in \text{RV}_\infty(0)$, we have 
\[
\lim_{t\to\infty} \frac{x(t)}{y(t)}
=\lim_{t\to\infty}\frac{x(t)}{F^{-1}(t)}=\lim_{t\to\infty}\frac{F^{-1}(\Lambda_\ast t)}{F^{-1}(t)} = 1. 
\]
so once again the solution of \eqref{eq.odepert} is of the same order as the solution of \eqref{eq.ode}. 

The proof of part (ii) of the theorem is identical in all respects to that of part (i), and therefore we present only the proof of part (i) in Section 4. In fact, there is a greater alignment of the hypotheses that appears at a first glance. When $f\in \text{RV}_0(\beta)$ for $\beta>1$, it follows that $F\in \text{RV}_0(1-\beta)$ 
and therefore that $F^{-1}\in \text{RV}_\infty(-1/(\beta-1))$ and $f\circ F^{-1}\in \text{RV}_\infty(-\beta/(\beta-1))$. Hence we see that the hypothesis of 
part (ii) are in some sense the limit of those in part (i) when $\beta\to \infty$. This suggests that part (ii) of the theorem applies in the case when $f$ is a rapidly varying function at $0$, and the solutions of the unperturbed differential equation are slowly varying at infinity. Moreover, the solution of 
the perturbed differential equation should also be slowly varying in this case. We present an example which supports these claims in the next section. First, 
we make some connections between the hypotheses in part (ii), especially with rapidly varying functions.  
\begin{remark}
Suppose $f\circ F^{-1}\in\text{RV}_\infty(-1)$. Then $F^{-1}\in \text{RV}_\infty(0)$. Therefore, we do not need to assume this second hypothesis in part (ii)
of Theorem~\ref{theorem.FxtLambda}.  
\end{remark}
\begin{proof}[Proof of Remark]
To see this, let $z'(t)=-f(z(t))$ for $t>0$ and $z(0)=1$. Then $z(t)=F^{-1}(t)$. Hence $0<-z'(t)=(f\circ F^{-1})(t)$. Therefore $-z'\in \text{RV}_\infty(-1)$, and 
so $z(t)-z(T)=\int_t^T -z'(s)\,ds$. Letting $T\to\infty$, we have $z(t)=\int_t^\infty -z'(s)\,ds$. Since $-z'\in \text{RV}_\infty(-1)$, it follows that 
$z\in\text{RV}_\infty(0)$. Hence $F^{-1}\in \text{RV}_\infty(0)$, as claimed. 
\end{proof}

\begin{remark}
Suppose that $f\in \text{RV}_0(\infty)$. Then $F^{-1}\in \text{RV}_\infty(0)$. 
\end{remark}
\begin{proof}[Proof of Remark]
The hypothesis that $f$ is rapidly varying at zero means by definition that
\[
\lim_{x\to 0^+} \frac{f(\lambda x)}{f(x)}=\left\{ 
\begin{array}{cc}
+\infty, & \lambda>1, \\
0, & \lambda<1.
\end{array}
\right.
\]
Now by the continuity of $f$ and l'H\^{o}pital's rule, we have 
\[
\lim_{x\to 0^+} \frac{F(\lambda x)}{F(x)}=\lim_{x\to 0^+}\frac{\int_{\lambda x}^1 \frac{1}{f(u)}\,du}{\int_{x}^1 \frac{1}{f(u)}\,du}
=\lim_{x\to 0^+} \frac{\lambda f(x)}{f(\lambda x)}
=
\left\{ 
\begin{array}{cc}
+\infty, & \lambda<1, \\
0, & \lambda>1.
\end{array}
\right.
\]
Consider the function $F_1(t)=F(1/t)$ as $t\to\infty$. Then 
\[
\lim_{t\to\infty} \frac{F_1(\lambda t)}{F_1(t)}=
\lim_{x\to 0^+} \frac{F(\lambda^{-1} x)}{F(x)}=
\left\{ 
\begin{array}{cc}
+\infty, & \lambda>1, \\
0, & \lambda<1.
\end{array}
\right.
\]
Therefore, $F_1$ is in $\text{RV}_\infty(\infty)$ and we have $\lim_{t\to\infty} F_1(t)=\lim_{x\to 0^+} F(x)=+\infty$, so 
$F_1^{-1}\in \text{RV}_\infty(0)$. Now $F^{-1}(F_1(t))=1/t$. Hence $F^{-1}(x)=F^{-1}(F_1(F_1^{-1}(x)))=1/F_1^{-1}(x)$. 
Therefore $F^{-1}\in \text{RV}_\infty(0)$.
\end{proof}

We notice that viewed as a function of $L$, $\Lambda_\ast:(0,\infty)\to (0,1)$ is decreasing and continuous with $\lim_{L\to 0^+} \Lambda_\ast(L)=1$ and 
$\lim_{L\to\infty} \Lambda_\ast(L)=0$. The first limit demonstrates that the limit in \eqref{eq.limFxtLambda} is a continuous extension of the limit observed in Theorem~\ref{theorem.gsmall}, because the hypothesis \eqref{eq.gfF0} can be viewed as \eqref{eq.gfFL} with $L=0$, while the 
resulting limiting behaviour of the solution  \eqref{eq.Fxtto1} can be viewed as \eqref{eq.limFxtLambda} where $\Lambda_\ast=1$. The monotonicity of $\Lambda_\ast$ in $L$ indicates that the slower the decay rate of the perturbation (i.e., the greater is $L$) the slower the rate of decay of the 
solution of \eqref{eq.odepert}. Since $\lim_{L\to\infty} \Lambda_\ast(L)=0$, this result also suggests that 
\begin{equation} \label{eq.gfFinfty}
\lim_{t\to\infty} \frac{g(t)}{(f\circ F^{-1})(t)}=+\infty
\end{equation}
implies 
\begin{equation} \label{eq.xslow}
\lim_{t\to\infty} \frac{F(x(t))}{t}=0,
\end{equation}
so that the solution of the perturbed differential equation entirely loses the decay properties of the underlying unperturbed equation when the perturbation 
$g$ exceeds the critical size indicated by \eqref{eq.gfFL}, and decays more slowly yet. This conjecture is borne out by virtue of the next theorem.
\begin{theorem}  \label{theorem.Fxtto0}
Suppose that $f$ obeys \eqref{eq.fglobalstable}, \eqref{eq.floclip} and that $F$ defined by \eqref{def.F} obeys \eqref{eq.Ftoinfty}.
Suppose further that $f$ and $g$ obey \eqref{eq.fgcns} and that \eqref{eq.gposxipos} holds. 
Let $x$ be the unique continuous solution of \eqref{eq.odepert}. 
Suppose that there exists $\phi$ such that \eqref{eq.fasyphimon} holds, and 
suppose further that $f$ and $g$ obey \eqref{eq.gfFinfty}. Suppose finally that $x(t)\to 0$ as $t\to\infty$. 
If $f\in \text{RV}_0(\beta)$ for $\beta>1$ or $f\circ F^{-1}\in \text{RV}_\infty(-1)$, then 
the unique continuous solution of \eqref{eq.odepert} obeys \eqref{eq.xslow}. 
\end{theorem}
\begin{remark}
We observe that the hypothesis that $x(t)\to 0$ as $t\to\infty$ has been appended to the theorem. This is because the slow rate of decay of $g$ may now 
cause solutions to tend to infinity, if coupled with a hypothesis on $f$ which forces $f(x)$ to tend to zero as $x\to\infty$ at a sufficiently rapid rate. We prefer to add this hypothesis, rather than sufficient conditions on $f$ and $g$ which would guarantee $x(t)\to\infty$.

We provide an example in which all the conditions of part (i) of Theorem~\ref{theorem.fxtgt1} are satisfied apart from $x(t)\to 0$ as $t\to\infty$, and show that it is in fact possible to get $x(t)\to\infty$. Let $\beta>1$ and $\theta\in (0,1)$ and consider the initial value problem 
\[
x'(t)=-f(x(t))+g(t), \quad t>0; \quad x(0)=1
\]
where $f(x)=x^\beta e^{-x}$ and $g(t)=(1-\theta)(1+t)^{-\theta}+(1+t)^{\beta(1-\theta)}e^{-(1+t)^{1-\theta}}$. 
The solution of this initial value problem is $x(t)=(1+t)^{1-\theta}$, so $x(t)\to\infty$ as $t\to\infty$. 
\end{remark}

We have shown, when \eqref{eq.gfFinfty} holds, that $F(x(t))/t\to 0$ as $t\to\infty$, so that the rate of decay of solutions of \eqref{eq.odepert} is slower 
than that of \eqref{eq.ode}. In the next theorem, under strengthened hypotheses on $g$, we determine the exact convergence rate to $0$ of the solution of \eqref{eq.odepert} when \eqref{eq.gfFinfty} holds, and we will show that the limit \eqref{eq.xslow} also holds. Once again, we add the hypothesis that $x(t)\to 0$
as $t\to\infty$.

\begin{theorem} \label{theorem.fxtgt1}
Suppose that $f$ obeys \eqref{eq.fglobalstable}, \eqref{eq.floclip} and that $F$ defined by \eqref{def.F} obeys \eqref{eq.Ftoinfty}.
Suppose further that $f$ and $g$ obey \eqref{eq.fgcns} and that \eqref{eq.gposxipos} holds. 

Suppose further that \eqref{eq.gfFinfty} holds 
and that $f\in\text{RV}_0(\beta)$ for some $\beta>1$ and $g\in \text{RV}_\infty(-\theta)$ for $\theta\geq  0$.
Let $x$ be the unique continuous solution of \eqref{eq.odepert} and suppose that $x(t)\to 0$ as $t\to\infty$. 
\begin{itemize}
\item[(i)] If $\theta> 0$, then
\begin{equation} \label{eq.fxtgtlim}
\lim_{t\to\infty} \frac{f(x(t))}{g(t)}=1.
\end{equation}
\item[(ii)] If $\theta=0$ and $g$ is asymptotic to a decreasing function, then $x$ obeys \eqref{eq.fxtgtlim}.
\end{itemize}
\end{theorem}
\begin{remark}
Unlike Theorem~\ref{theorem.fxtgt1}, previous theorems have not assumed that $g$ is regularly varying, or obeys other regular asymptotic properties, beyond asking that $g$ decays in some manner related to $f\circ F^{-1}$. However, the assumption that $g$ is regularly (or slowly varying) in Theorem~\ref{theorem.fxtgt1} is quite natural, as by \eqref{eq.gfFinfty} it decays more slowly to zero than a function which is itself regularly varying at infinity (with negative index $-\beta/(\beta-1)$). Moreover, it is a consequence of the hypotheses of Theorem~\ref{theorem.FxtLambda} that $g$ is regularly varying, as it is asymptotic to 
$f\circ F^{-1}$ which is assumed to be regularly varying (with index $-\beta/(\beta-1)$ in part (i), and index -1 in part (ii)). We notice moreover that 
Theorem~\ref{theorem.fxtgt1} does not deal with the case when $f$ is rapidly varying at $0$.
\end{remark}

\begin{remark}
It is interesting to note that \eqref{eq.fxtgtlim} may be thought of as \eqref{eq.limFxtLambda} in the limit $L\to\infty$. To see this, notice if 
\eqref{eq.limFxtLambda} holds, we have 
\[
\lim_{t\to\infty} \frac{g(t)}{f(F^{-1}(t))}=L, \quad \lim_{t\to\infty} \frac{F(x(t))}{\Lambda_\ast(L)t}=1.
\]
Therefore, if $f\circ F^{-1}$ is regularly varying, we have
\[
\lim_{t\to\infty} \frac{f(x(t))}{(f\circ F^{-1})(\Lambda_\ast t)}=
\lim_{t\to\infty} \frac{(f\circ F^{-1})(F(x(t))}{(f\circ F^{-1})(\Lambda_\ast t)}=1.
\]
Therefore, if $\beta>1$, we have
\[
\lim_{t\to\infty} \frac{(f\circ F^{-1})(\Lambda_\ast t)}{(f\circ F^{-1})(t)}=\Lambda_\ast^{-\beta/(\beta-1)}.
\]
Hence
\[
\lim_{t\to\infty} \frac{f(x(t))}{(f\circ F^{-1})(t)}= 
\lim_{t\to\infty} \frac{f(x(t))}{(f\circ F^{-1})(\Lambda_\ast t)}\cdot \frac{(f\circ F^{-1})(\Lambda_\ast t)}{(f\circ F^{-1})(t)}
=\Lambda_\ast^{-\beta/(\beta-1)}.
\]
Therefore
\begin{equation} \label{eq.Ltoinfty}
\lim_{t\to\infty} \frac{f(x(t))}{g(t)}=\lim_{t\to\infty} \frac{f(x(t))}{(f\circ F^{-1})(t)}\cdot \frac{(f\circ F^{-1})(t)}{g(t)}=\frac{1}{L} \Lambda_\ast^{-\beta/(\beta-1)}.
\end{equation}
In case (i) of Theorem~\ref{theorem.FxtLambda}, we have $\beta>1$ and 
$(1-\Lambda_\ast)\Lambda^{-\beta/(\beta-1)}=L$, so $1-\Lambda_\ast=L/\Lambda_\ast^{-\beta/(\beta-1)}$. Since $\Lambda_\ast(L)\to 0$ as 
as $L\to\infty$, we have 
\[
\lim_{L\to\infty} \frac{L}{\Lambda_\ast^{-\beta/(\beta-1)}}=1,
\]
and therefore the limit on the righthand side of \eqref{eq.Ltoinfty} as $L\to \infty$ is unity.
In case (ii) of Theorem~\ref{theorem.FxtLambda}, in place of \eqref{eq.Ltoinfty} we find that 
\[
\lim_{t\to\infty} \frac{f(x(t))}{g(t)}=\frac{1}{L} \Lambda_\ast^{-1}.
\]
Since $\Lambda_\ast(L)=1/(1+L)$, we have that the righthand side once again tends to unity as $L\to\infty$. Therefore, we see that 
the rate of decay changes smoothly as the parameter $L$ changes from being zero, to finite, and then to infinity. 
\end{remark}

\begin{remark}
We remark that under the hypotheses of Theorem~\ref{theorem.fxtgt1}, we have the limit $F(x(t))/t\to 0$ as $t\to\infty$, which is consistent with the result of Theorem~\ref{theorem.Fxtto0}. To see this, we note that 
under the hypothesis \eqref{eq.gfFinfty}, we conclude that \eqref{eq.fxtgtlim}. Multiplying these limits gives   
\[
\lim_{t\to\infty} \frac{f(x(t))}{f(F^{-1}(t))}=+\infty.
\]
Since $\beta>1$, using the fact that $f\in\text{RV}_0(\beta)$ we have 
$\lim_{t\to\infty} x(t)/F^{-1}(t)=+\infty$, and since $F\in \text{RV}_\infty(1-\beta)$ with $\beta>1$ and $F$ is decreasing, we obtain the limit  
$\lim_{t\to\infty} F(x(t))/t=0$, as required. 
\end{remark}
We may now consolidate our findings into two theorems which characterise the asymptotic behaviour of solutions of \eqref{eq.odepert}: in the first, we make no assumption about the regular or slow variation of $g$ at infinity, and allow $f$ to be regularly or rapidly varying at zero; in the second, we assume that both $f$ and $g$ are both regularly varying, and obtain exact asymptotic estimates on the solution in each case. 
\begin{theorem} \label{theorem.overallnorv}
Suppose that $f$ obeys \eqref{eq.fglobalstable}, \eqref{eq.floclip} and that $F$ is defined by \eqref{def.F}.
Suppose also that $f$ and $g$ obey \eqref{eq.fgcns} and that \eqref{eq.gposxipos} holds. 
Suppose that $f\in\text{RV}_0(\beta)$ for some $\beta>1$ or that $f\circ F^{-1}\in \text{RV}_\infty(-1)$, and that 
$f$ obeys \eqref{eq.fasyphimon}. Let $x$ be the unique continuous solution of \eqref{eq.odepert} and suppose that $x(t)\to 0$ as $t\to\infty$. Finally, 
suppose that
\[
\lim_{t\to\infty} \frac{g(t)}{f(F^{-1}(t))}=L\in [0,\infty].
\]
\begin{itemize}
\item[(i)] If $L=0$, then 
\[
\lim_{t\to\infty} \frac{F(x(t))}{t}=1.
\]
\item[(ii)] If $L\in (0,\infty)$, then 
\[
\lim_{t\to\infty} \frac{F(x(t))}{t}=\Lambda_\ast(L),
\]
where $\Lambda_\ast\in (0,1)$ is given by 
(I) the unique solution of $(1-\Lambda_\ast)\Lambda_\ast^{-\beta/(\beta-1)}=L$ when $f\in \text{RV}_0(\beta)$ for some $\beta>1$ 
and (II) $\Lambda^\ast=1/(1+L)$ when $f\circ F^{-1}\in \text{RV}_\infty(-1)$.
\item[(iii)] If $L=\infty$, then 
\[
\lim_{t\to\infty} \frac{F(x(t))}{t}=0.
\]
\end{itemize}
\end{theorem}
Theorem~\ref{theorem.overallnorv} is established by combining the results of Theorems~\ref{theorem.gsmall},~\ref{theorem.FxtLambda} and~\ref{theorem.Fxtto0}.
On the other hand, by combining the results of Theorems~\ref{theorem.gsmall},~\ref{theorem.FxtLambda} and~\ref{theorem.fxtgt1}, we arrive at a classification of the dynamics of \eqref{eq.odepert} when $f$ and $g$ are regularly varying.
\begin{theorem}  \label{theorem.overall}
Suppose that $f$ obeys \eqref{eq.fglobalstable}, \eqref{eq.floclip} and that $F$ is defined by \eqref{def.F}.
Suppose also that $f$ and $g$ obey \eqref{eq.fgcns} and that \eqref{eq.gposxipos} holds. 
Suppose that $f\in\text{RV}_0(\beta)$ for some $\beta>1$ and that $g\in \text{RV}_\infty(-\theta)$ for $\theta> 0$.
Let $x$ be the unique continuous solution of \eqref{eq.odepert} and suppose that $x(t)\to 0$ as $t\to\infty$. Finally, 
suppose that
\[
\lim_{t\to\infty} \frac{g(t)}{f(F^{-1}(t))}=L\in [0,\infty].
\]
\begin{itemize}
\item[(i)] If $L=0$, then 
\[
\lim_{t\to\infty} \frac{F(x(t))}{t}=1.
\]
\item[(ii)] If $L\in (0,\infty)$, then 
\[
\lim_{t\to\infty} \frac{F(x(t))}{t}=\Lambda_\ast(L),
\]
where $\Lambda_\ast\in (0,1)$ is the unique solution of $(1-\Lambda_\ast)\Lambda_\ast^{-\beta/(\beta-1)}=L$.
\item[(iii)] If $L=\infty$, then 
\[
\lim_{t\to\infty} \frac{f(x(t))}{g(t)}=1.
\]
\end{itemize}
\end{theorem}
We close by remarking that in cases (i) and (ii), the solution of \eqref{eq.odepert} is regularly varying at infinity with index $-\beta/(\beta-1)$, 
while in case (iii) it is regularly varying at infinity with index $-\theta/\beta$. 

\section{Examples}
We next demonstrate the scope of the theorems by studying a number of examples. We have expressly chosen the examples so that solutions are known in closed form. This enables us to demonstrate independently of our theorems the breadth of the results in the paper.
 
We start with an example that demonstrates that when $g(t)$ does not have the same sign as the initial condition $\xi$, and the solution $x$ of \eqref{eq.odepert} nonetheless retains the sign of the initial condition, the perturbation $g$ can be small in the sense that \eqref{eq.gfF0} holds, but the solution $x$ of \eqref{eq.odepert} does not obey \eqref{eq.Fxtto1}. This shows the importance of retaining the assumption that $g$ be positive in Theorem~\ref{theorem.gsmall}.
\begin{example} \label{examp.gneg}
Suppose that $\beta>1$, $\theta>\beta/(\beta-1)$. Let $\xi\in (0,(\theta-1)^{1/(\beta-1)})$. Suppose that 
\[
g(t)=-(1+t)^{-\theta}\left(\xi(\theta-1)-\xi^\beta(1+t)^{-\beta(\theta-1)+\theta}\right), \quad t\geq 0.
\]
Notice that $g(t)<0$ for all $t\geq 0$. Let $f(x)=\sgn(x)|x|^\beta$ for $x\in\mathbb{R}$. 
Then the unique continuous solution of \eqref{eq.odepert} is $x(t)=\xi(1+t)^{-(\theta-1)}$ for $t\geq 0$. In the terminology of this paper, we have 
\[
F(x)= \frac{1}{\beta-1}\left(x^{-\beta+1}-1\right),
\] 
so $\lim_{x\to 0^+} F(x)/x^{-\beta+1}=1/(\beta-1)$. 
Hence $F^{-1}(t) \sim \left((\beta-1)t\right)^{-1/(\beta-1)}$ as $t\to\infty$ and so 
$(f\circ F^{-1})(t)\sim \left((\beta-1)t\right)^{-\beta/(\beta-1)}$
as $t\to\infty$. Since $\theta>\beta/(\beta-1)$, it follows that $g$ and $f$ obeys \eqref{eq.gfF0}. However, 
\[
\lim_{t\to\infty} \frac{F(x(t))}{t}= \frac{1}{\beta-1}\lim_{t\to\infty} \frac{x(t)^{-\beta+1}}{t}
= \frac{\xi^{-\beta+1}}{\beta-1}\lim_{t\to\infty} t^{(\theta-1)(\beta-1)-1}  =+\infty,
\]
so the conclusion of Theorem~\ref{theorem.gsmall} does not hold.
\end{example}
The next example concerns an equation of the form \eqref{eq.odepert} to which Theorem~\ref{theorem.gsmall} could be applied, but for which a closed form solution is known, and therefore independently exemplifies this theorem.
\begin{example} \label{examp.g0}
Let $\eta>0$, $\beta>1$, $\xi>0$, and let $A=\xi^{1-\beta}$. Suppose that 
\[
g(t)=\frac{\eta A(1+t)^{-(\eta+1)}}{\beta-1}\left\{A(1+t)^{-\eta}+(\beta-1)t \right\}^{-\beta/(\beta-1)}, \quad t\geq 0.
\]
Then $g(t)>0$ for all $t\geq 0$. Suppose that $f(x)=\sgn(x)|x|^\beta$ for $x\geq 0$. Then the unique continuous solution of the initial value 
problem \eqref{eq.odepert} is 
\[
x(t)=\left( A(1+t)^{-\eta}  + (\beta-1) t  \right)^{-1/(\beta-1)}
\]
We notice that   
\[
\lim_{t\to\infty} \frac{g(t)}{t^{-(\eta+1+\beta/(\beta-1))}}=\frac{\eta A}{\beta-1} (\beta-1)^{-\beta/(\beta-1)}.
\]
so, as $\eta>0$, we have that $g$ obeys \eqref{eq.gfF0}. It can be seen that all the hypotheses of Theorem~\ref{theorem.gsmall} hold. 
On the other hand, from the definition of $F$ we have that \eqref{eq.Fxtto1} holds which we are able to conclude independently of Theorem~\ref{theorem.gsmall}.  
\end{example}
We now give an example to which part (i) of Theorem~\ref{theorem.FxtLambda} applies.
\begin{example} \label{examp.Lgt0}
Suppose that $A>1/(\beta-1)^{1/(\beta-1)}$. Define $\xi>0$ and 
\[
g(t)=\left(A^\beta-A\frac{1}{\beta-1}\right)\left\{\left(\frac{A}{\xi}\right)^{\beta-1}+t\right\}^{-\beta/(\beta-1)}, \quad t\geq 0.
\] 
Then $g(t)>0$ for all $t\geq 0$. Suppose also that $f(x)=\sgn(x)|x|^\beta$ for $x\geq 0$. Then the initial value problem \eqref{eq.odepert} has unique continuous solution $x(t)=A( (A/\xi)^{\beta-1}+t)^{-1/(\beta-1)}$ for $t\geq 0$. Notice that 
\[
\lim_{t\to\infty} \frac{g(t)}{(f\circ F^{-1})(t)}=\frac{A^\beta-\frac{1}{\beta-1}A}{(\beta-1)^{-\beta/(\beta-1)}}=:L>0.
\]
Also we have that 
\[
\lim_{t\to\infty} \frac{x(t)}{F^{-1}(t)}=\frac{A}{(\beta-1)^{-1/(\beta-1)}},
\]
so 
\[
\lim_{t\to\infty} \frac{F(x(t))}{t}=\frac{A^{1-\beta}}{\beta-1}=:\Lambda_\ast.
\]
Since $A>1/(\beta-1)^{1/(\beta-1)}$, we have $\Lambda_\ast\in (0,1)$ and moreover one can check that $(1-\Lambda_\ast)\Lambda_\ast^{-\beta/(\beta-1)}=L$.
Therefore it can be seen that the conclusion of Theorem~\ref{theorem.FxtLambda} applies. 
\end{example}
Even though our results cover more comprehensively the case when $f$ has ``power--like'' behaviour close to zero, our next example 
demonstrates that when $f$ is rapidly varying at zero (and in fact has all its one--sided derivatives equal to 0 at 0), we can still determine the 
rate of convergence of solutions. Theorem~\ref{theorem.FxtLambda} part (ii) covers this example. 
\begin{example} \label{examp.Lgt0frapid}
Suppose that $f(x)=\sgn(x)e^{-1/|x|}$ for $x\neq 0$ and $f(0)=0$. Then for $x>0$ we have
\[
F(x)=\int_x^1 e^{1/u}\,du = \int_{1}^{1/x} v^{-2}e^v\,dv.
\]
Therefore by l'H\^opital's rule we have
\[
\lim_{x\to 0^+} \frac{F(x)}{e^{1/x}x^2} = \lim_{y\to \infty} \frac{\int_1^y v^{-2}e^v\,dv}{e^y y^{-2}} 
= \lim_{y\to\infty} \frac{e^y y^{-2}}{e^y y^{-2}-2y^{-3}e^y}=1.
\]
Since $F^{-1}(t)\to 0$ as $t\to\infty$ we have
\begin{equation} \label{eq.e1xFx}
\lim_{t\to\infty} \frac{t}{e^{1/F^{-1}(t)}F^{-1}(t)^2}=1.
\end{equation}
Therefore 
\[
\lim_{t\to\infty} \left\{\log t - \frac{1}{F^{-1}(t)} - 2 \log F^{-1}(t)\right\} = 0.
\]
Since $\lim_{x\to 0^+} \log (x)/x^{-1} = \lim_{x\to 0^+} x \log (x) = 0$, we have 
\[
\lim_{t\to\infty} \frac{\log t}{\frac{1}{F^{-1}(t)}}=1,
\]
so $F^{-1}(t)/(\log t)^{-1}\to 1$ as $t\to\infty$. Hence $F^{-1}\in \text{RV}_\infty(0)$. Moreover, 
since $f(F^{-1}(t))=e^{-1/F^{-1}(t)}$, from \eqref{eq.e1xFx} we obtain 
\[
\lim_{t\to\infty} \frac{(f\circ F^{-1})(t)}{F^{-1}(t)^2/t}=1,
\]
and therefore, as $F^{-1}(t)/(\log t)^{-1}\to 1$ as $t\to\infty$, we get 
\[
\lim_{t\to\infty} \frac{(f\circ F^{-1})(t)}{\frac{1}{t(\log t)^2}}=1.
\]
Hence $f\circ F^{-1}\in \text{RV}_\infty(-1)$.

Define
\begin{multline} \label{eq.gfrapid}
\tilde{g}(t)=\frac{3}{(3e+t)\log^2(e+t)}-\frac{1}{(3e+t)\log^2((e+t/3)\log^2(e+t))}
\\-\frac{2}{(e+t)\log(e+t)\log^2((e+t/3)\log^2(e+t))}, \quad t\geq 3.
\end{multline}
Notice that $\tilde{g}$ is continuous and positive on $[3,\infty)$. Then we have that 
\[
\tilde{x}(t)=\frac{1}{\log\left((e+t/3)\log^2(e+t)\right)}, \quad t\geq 3
\]
is a solution of the initial value problem 
\begin{equation} \label{eq.xtildeode}
\tilde{x}'(t)=-f(\tilde{x}(t))+\tilde{g}(t), \quad t>3; \quad \tilde{x}(3)=\xi=\frac{1}{\log(e+1)\log^2(e+3)}>0.
\end{equation}
Now define $x(t)=\tilde{x}(t+3)$ and $g(t)=\tilde{g}(t+3)$ for $t\geq 0$. Then $g$ is continuous and positive on $[0,\infty)$
and $x$ satisfies the initial value problem
\[
x'(t)=-f(x(t))+g(t), \quad t>0; \quad x(0)=\xi=\frac{1}{\log(e+1)\log^2(e+3)}>0.
\]

To see that $\tilde{x}$ obeys \eqref{eq.xtildeode}, define
\[
\eta(t)=\frac{\log(e+t)}{\log\left((e+t/3)\log^2(e+t)\right)}, \quad t\geq 3.
\]
Then $e^{1/\tilde{x}(t)}=(e+t)^{1/\eta(t)}$ so by the definition of $\tilde{x}$ we have 
\[
(e+t)^{1/\eta(t)}=\frac{1}{3}(3e+t)\log^2(e+t), \quad t\geq 3.
\]
Also 
$f(\tilde{x}(t))=1/(e+t)^{1/\eta(t)}=3/((3e+t)\log^2(e+t))$. This is the first term on the righthand side of 
\eqref{eq.gfrapid}. It is easy to check directly from the formula for $\tilde{x}$ that the second and third terms on the righthand side equal $\tilde{x}'(t)$. 
Therefore $\tilde{g}(t)=f(\tilde{x}(t))+\tilde{x}'(t)$, so $\tilde{x}$ obeys \eqref{eq.xtildeode}.

Notice that 
\[
\lim_{t\to\infty} \frac{g(t)}{(f\circ F^{-1})(t)}=2.
\]
We can determine the asymptotic behaviour of $x(t)$ as $t\to\infty$ using the auxiliary function $\eta$. 
Since 
\[
e^{1/x(t)}=e^{1/\tilde{x}(t+3)}
=
(e+t+3)^{1/\eta(t+3)}=\frac{1}{3}(3e+t+3)\log^2(e+t+3),
\]
and $x(t)/(\log t)^{-1}\to 1$ as $t\to\infty$, we can check that 
\[
\lim_{t\to\infty} \frac{F(x(t))}{t}
=\lim_{t\to\infty} \frac{e^{1/x(t)}x(t)^2}{t}
=\lim_{t\to\infty} \frac{\frac{1}{3}(3e+t+3)\log^2(e+t+3) (\log t)^{-2}}{t}
=\frac{1}{3}.
\]
This calculation is independent of Theorem~\ref{theorem.FxtLambda} part (ii) but confirms it, because here $L=2$ and $\Lambda_\ast=1/(L+1)=1/3$.
\end{example}

We now present an example to which Theorem~\ref{theorem.fxtgt1} applies.
\begin{example} \label{examp.Linfty}
Let $\theta<\beta/(\beta-1)$ and $\xi<(\beta/\theta)^{1/(1-\beta+\beta/\theta)}$. Suppose that 
\[
g(t)=(\xi^{-\beta/\theta}+t)^{-\theta}\left(1-\frac{\theta}{\beta}(\xi^{-\beta/\theta}+t)^{-\frac{\theta}{\beta}-1+\theta}\right), \quad t\geq 0.
\]
Notice that $g(t)>0$ for all $t\geq 0$. Let $f(x)=\sgn(x)|x|^\beta$ for $x\geq 0$. 
Then the unique solution of the initial value problem \eqref{eq.odepert} is $x(t)=(\xi^{-\beta/\theta}+t)^{-\theta/\beta}$ for 
$t\geq 0$. We see that $g\in\text{RV}_\infty(-\theta)$, and also that $g$ and $f$ obey 
\eqref{eq.gfFinfty}. Also $x(t)\to 0$ as $t\to\infty$. Hence all the hypotheses of Theorem~\ref{theorem.fxtgt1} hold. Moreover, we can see, 
independently of the conclusion of Theorem~\ref{theorem.fxtgt1} that \eqref{eq.fxtgtlim} holds.
\end{example}
Our final example shows how Theorem~\ref{theorem.fxtgt1} can detect very slowly decaying (i.e., slowly varying) solutions of \eqref{eq.odepert}. 
This arises when the perturbation $g$ is slowly varying at infinity. In our example, the perturbation exhibits iterated logarithmic decay. 
\begin{example}
Let $\xi>(\beta e^{e+1})^{-1/(\beta-1)}$. Define $f(x)=|x|^\beta \sgn(x)$ for $x\geq 0$ and 
\[
g(t)=\frac{\xi^\beta}{\log_2(t+e^e)}- \xi\frac{1}{\beta}\cdot \frac{1}{(\log_2(t+e^e))^{1+1/\beta}}\cdot \frac{1}{(t+e^e)\log(t+e^e)}, \quad t\geq 0.
\] 
Notice that the restriction on $\xi$ implies that $g(t)>0$ for all $t\geq 0$. It can then be verified that   
\[
x(t)=\frac{\xi}{(\log_2(t+e^e))^{1/\beta}}, \quad t \geq 0
\]
is the unique continuous solution of \eqref{eq.odepert}. Notice that 
\[
\lim_{t\to\infty} \frac{g(t)}{\frac{\xi^\beta}{\log_2(t+e^e)}}=1
\]
and that 
\[
\lim_{t\to\infty} \frac{f(x(t))}{g(t)}= \lim_{t\to\infty} \frac{\xi^\beta (\log_2(t+e^e))^{-1} }{\xi^\beta (\log_2(t+e^e)^{-1}}=1.
\]
Therefore once again, we can confirm the conclusion of Theorem~\ref{theorem.fxtgt1} independently of its proof. Of course, it can be shown that all the hypotheses of part (ii) of Theorem~\ref{theorem.fxtgt1} hold for this problem; in particular, we may take the decreasing function to which $g$ is asymptotic to be to be  $\gamma(t)=\xi^\beta/\log_2(t+e^e)$ for $t\geq 0$. 
\end{example}

\section{Proofs}
\subsection{Proof of Theorem~\ref{theorem.gsmall}}
Since $F^{-1}(0)=1$, we have that $F^{-1}(t)\in (0,1)$ for all $t> 0$. Let $T>0$. Then
\[
\int_0^T f(F^{-1}(s))\,ds =  \int_{F^{-1}(0)}^{F^{-1}(T)}f(u)\cdot \left(-\frac{1}{f(u)}\right)\,du = 1-F^{-1}(T).
\]
Since $F^{-1}(t)\to 0$ as $t\to\infty$, we have that $f\circ F^{-1}\in L^1(0,\infty)$. Therefore as \eqref{eq.gfF0} holds we have
that $g\in L^1(0,\infty)$. Therefore, by virtue of Proposition~\ref{prop.xto0ginL1} we have that $x(t)\to 0$ as $t\to\infty$.

Suppose that $\xi>1$. Since $x(t)\to 0$ as $t\to\infty$, it follows that there exists $0<T_1:=\sup\{t>0:x(t)=1\}$. Define $z$ by
$z'(t)=-f(z(t))$ for $t>T_1$ and $z(T_1)=1$. Then $x'(T_1)=-f(z(T_1))+g(T_1)>z'(T_1)$, because $g(t)>0$ for all $t>0$. Therefore,
we have that $z(t)<x(t)$ for all $t>T_1$. If $\delta\geq 1$, notice that $\phi(x(t))>\phi(z(t))$ for all $t>T_1$.
If $\delta<1$, we notice that there is $T_1<T_2:=\sup\{t>0:x(t)=\delta\}$. Moreover, $z(t)<\delta$ for $t>T_2$. Hence
for $t>T_2$ we have $\phi(x(t))>\phi(z(t))$ and $z(t)<\delta$ for all $t>T_2$.
Therefore, irrespective of the level of $\delta$, there exists $T_3> T_1$ such that $\phi(x(t))>\phi(z(t))$ and $z(t)<\delta$ for all $t\geq T_3$.
Next, $z$ is given by $z(t)=F^{-1}(t-T_1)$ for all $t\geq T_3\geq T_1$. Also $t\mapsto (\phi\circ F^{-1})(t-T_1)$ is decreasing on $[T_3,\infty)$, and
so for $t>T_1+T_3$ we must have that $(\phi\circ F^{-1})(t-T_1)>(\phi\circ F^{-1})(t)$. Therefore, for $t>T_1+T_3$ it follows that
\[
\frac{g(t)}{\phi(x(t))}<\frac{g(t)}{\phi(z(t))}=\frac{g(t)}{(\phi\circ F^{-1})(t-T_1)}<\frac{g(t)}{(\phi\circ F^{-1})(t)}.
\]
Therefore by \eqref{eq.gfF0}, we have that
\[
\lim_{t\to\infty} \frac{g(t)}{\phi(x(t))}=0.
\]
Since $x(t)\to 0$ as $t\to\infty$ and $f(x)/\phi(x)\to 1$ as $x\to 0^+$, we have $g(t)/f(x(t))\to 0$ as $t\to\infty$. Hence
\[
\lim_{t\to\infty} \frac{x'(t)}{f(x(t))}=-1.
\]
Integration yields
\[
\lim_{t\to\infty} \frac{F(x(t))}{t}=1.
\]

In the case that $\xi\leq 1$, define $u$ to be the unique continuous solution of
\[
u'(t)=-f(u(t))+2g(t), \quad t>0; \quad u(0)=\xi+1.
\]
Then it can be shown by contradiction that $x(t)<u(t)$ for all $t\geq 0$.
Since $F$ is decreasing, we have $F(x(t))>F(u(t))$ for all $t\geq 0$. We may apply the argument for $\xi>1$ 
above to show that $F(u(t))/t\to 1$ as $t\to\infty$. Therefore, we have that
\[
\liminf_{t\to\infty} \frac{F(x(t))}{t}\geq 1.
\]
On the other hand, define $z$ by $z'(t)=-f(z(t))$ for $t>0$ and $z(0):=\xi'<\xi$. Then $x(t)>z(t)$ for $t\geq 0$. Therefore we have $F(x(t))\leq F(z(t))$ for
all $t\geq 0$. However, $F(z(t))-F(\xi')=t$, so $F(z(t))/t\to 1$ as $t\to\infty$. Therefore we have
\[
\limsup_{t\to\infty} \frac{F(x(t))}{t}\leq 1.
\]
Combining this with the limit inferior gives $F(x(t))/t\to 1$ as $t\to\infty$ as required.

\subsection{Proof of Theorem~\ref{theorem.FxtLambda}}
Let $\epsilon\in (0,1)$ and define $h(x)=x^{-\beta/(\beta-1)}(1-x)$ for $x\in (0,1]$. Then $h(x)\uparrow \infty$ as $x\downarrow 0$ and $h(1)=0$, with $h$
decreasing on $(0,1]$. Therefore, for each $L>0$ there exists a unique $\Lambda(\epsilon)\in (0,1)$ such that
\[
h(\Lambda(\epsilon))=L\frac{1+\epsilon}{(1-\epsilon)^2}.
\]
Since $h(\Lambda_\ast)=L$, the continuity of $h$ ensures that $\Lambda(\epsilon)\to \Lambda_\ast$ as $\epsilon\to 0^+$.

We have already shown in the proof of Theorem~\ref{theorem.gsmall} that $f\circ F^{-1}\in L^1(0,\infty)$. Since $g(t)/(f\circ F^{-1})(t)\to L$ as $t\to\infty$, it follows that $g\in L^1(0,\infty)$.
Therefore it follows from Proposition~\ref{prop.xto0ginL1} that $x(t)\to 0$ as $t\to\infty$.

Now, there is a $\delta<1$ and an increasing $\phi:(0,\delta)\to (0,\infty)$ such that $\phi(x)/f(x)\to 1$ as $x\to 0^+$. 
Hence for every $\epsilon\in (0,1)$ there exists $x_2(\epsilon)>0$ such that 
\[
\frac{f(x)}{\phi(x)}>1-\epsilon \quad\text{for all $x\in (0,x_2(\epsilon))$}.
\]
Since $F^{-1}(t)\to\infty$ as $t\to\infty$, there exists $T_1(\epsilon)>0$ such that
\[
g(t)<L(1+\epsilon)(\phi\circ F^{-1})(t), \quad t\geq T_1(\epsilon).
\]
There also exists $T_2(\epsilon)>0$ such that for $t\geq T_2(\epsilon)$ we have $x(t)< 1$ and $x(t)<x_2(\epsilon)$. Let $T_3(\epsilon)=\max(T_1(\epsilon),T_2(\epsilon))$. Since $f\circ F^{-1}\in \text{RV}_\infty(-\beta/(\beta-1))$, we have that
$\phi\circ F^{-1}\in \text{RV}_\infty(-\beta/(\beta-1))$. Therefore, for each $\eta\in (0,\epsilon)\subset (0,1)$ there exists $x_1'(\eta)>0$
such that
\[
\frac{(\phi\circ F^{-1})(\Lambda(\epsilon)x)}{(\phi\circ F^{-1})(x)}\geq \Lambda(\epsilon)^{-\beta/(\beta-1)}(1-\eta), \quad x\geq x_1'(\eta).
\]
Now fix $\eta=\epsilon/2$. Then
\[
\frac{(\phi\circ F^{-1})(\Lambda(\epsilon)x)}{(\phi\circ F^{-1})(x)}\geq \Lambda(\epsilon)^{-\beta/(\beta-1)}(1-\epsilon/2), \quad x\geq x_1'(\epsilon/2).
\]
If $x_1'(\epsilon/2)\leq F(\delta)$, define $x_1''(\epsilon)=F(\delta)+1$, so $x_1''(\epsilon)>F(\delta)$ and
\begin{equation} \label{eq.x1}
\frac{(\phi\circ F^{-1})(\Lambda(\epsilon)x)}{(\phi\circ F^{-1})(x)}\geq \Lambda(\epsilon)^{-\beta/(\beta-1)}(1-\epsilon/2), \quad x\geq x_1''(\epsilon).
\end{equation}
If $x_1'(\epsilon/2)>F(\delta)$, define $x_1''(\epsilon)=x_1'(\epsilon/2)$, so once again $x_1''(\epsilon)>F(\delta)$ and once again \eqref{eq.x1} holds.
If $x_1''(\epsilon)>F(x_2(\epsilon))/\Lambda(\epsilon)$, define $x_1(\epsilon)=x_1''(\epsilon)>F(\delta)$, so that
\begin{equation} \label{eq.x11}
\frac{(\phi\circ F^{-1})(\Lambda(\epsilon)x)}{(\phi\circ F^{-1})(x)}\geq \Lambda(\epsilon)^{-\beta/(\beta-1)}(1-\epsilon/2), \quad x\geq x_1(\epsilon).
\end{equation}
and
\begin{equation} \label{eq.f1x1x2}
F^{-1}(\Lambda(\epsilon) x_1(\epsilon))<x_2(\epsilon).
\end{equation}
On the other hand, if $x_1''(\epsilon)<F(x_2(\epsilon))/\Lambda(\epsilon)$, define $x_1(\epsilon)=F(x_2(\epsilon))/\Lambda(\epsilon)+1>F(x_2(\epsilon))/\Lambda(\epsilon)$, so \eqref{eq.f1x1x2} and \eqref{eq.x11} hold.
Moreover,
\begin{equation} \label{eq.x1Fdel}
x_1(\epsilon)>x_1''(\epsilon)>F(\delta).
\end{equation}

Define next
\[
T_4(\epsilon)=\sup\left\{t>0: x(t)=\frac{1}{2}F^{-1}(\Lambda(\epsilon)x_1(\epsilon))\right\}.
\]
If $x(t)<\frac{1}{2}F^{-1}(\Lambda(\epsilon)x_1(\epsilon))$
for all $t\geq 0$, set $T_4(\epsilon)=T_3(\epsilon)$. Finally, set $T(\epsilon)=1+\max(T_3(\epsilon),T_4(\epsilon),x_1(\epsilon))$. Since $T(\epsilon)>T_4(\epsilon)$, we have
\[
x(T(\epsilon))\leq \frac{1}{2}F^{-1}(\Lambda(\epsilon)x_1(\epsilon)),
\]
so $2x(T(\epsilon))\leq F^{-1}(\Lambda(\epsilon)x_1(\epsilon))$. 
Define
\[
M=\frac{F^{-1}(\Lambda(\epsilon)x_1(\epsilon))}
{x(T(\epsilon))}.
\]
Then $M\geq 2>1$. Finally, define
\[
x_U(t)=F^{-1}\biggl(\Lambda(\epsilon)(t-T(\epsilon))+F(Mx(T(\epsilon)))\biggr), \quad t\geq T(\epsilon).
\]
Thus $x_U(T(\epsilon))=Mx(T(\epsilon))>x(T(\epsilon))$. 

Since $F^{-1}(\Lambda(\epsilon) x_1(\epsilon))<x_2(\epsilon)$, we have $Mx(T(\epsilon))=F^{-1}(\Lambda(\epsilon) x_1(\epsilon))<x_2(\epsilon)$.
Hence $x_U(T(\epsilon))=Mx(T(\epsilon))<x_2(\epsilon)$. Therefore, it is always the case that $x_U(t)<x_2(\epsilon)$ for all $t\geq T(\epsilon)$.

Next, for $t>T(\epsilon)$, we have $F'(x_U(t))x_U'(t)=\Lambda(\epsilon)$. Therefore for $t>T(\epsilon)$ we have $x_U'(t)=-\Lambda(\epsilon)f(x_U(t))$, or
\[
x_U'(t)=-f(x_U(t))+(1-\Lambda(\epsilon))f(x_U(t))>-f(x_U(t)) +  (1-\Lambda(\epsilon))(1-\epsilon)\phi(x_U(t)).
\]
Since $M x(T(\epsilon))= F^{-1}(\Lambda(\epsilon)x_1(\epsilon))$, we get
$F(M x(T(\epsilon)))/\Lambda(\epsilon)=x_1(\epsilon)$. Therefore for $t\geq T(\epsilon)$, we have
\[
t-T(\epsilon)+\frac{F(Mx(T(\epsilon)))}{\Lambda(\epsilon)}\geq x_1(\epsilon).
\]
Hence for $t>T(\epsilon)$, by using this inequality and \eqref{eq.x11}, we have
\begin{align*}
\lefteqn{(1-\Lambda(\epsilon))(1-\epsilon)\phi(x_U(t))}\\
&=
(1-\Lambda(\epsilon))(1-\epsilon) (\phi \circ F^{-1})\left(\Lambda(\epsilon)\left\{t-T(\epsilon)+\frac{F(Mx(T(\epsilon)))}{\Lambda(\epsilon)}\right\}\right)\\
&\geq (1-\Lambda(\epsilon))(1-\epsilon)(1-\epsilon/2)\Lambda(\epsilon)^{-\beta/(\beta-1)} (\phi\circ F^{-1})\left(t-T(\epsilon)+\frac{F(Mx(T(\epsilon)))}{\Lambda(\epsilon)}  \right)\\
&=h(\Lambda(\epsilon))(1-\epsilon/2)(1-\epsilon)  (\phi\circ F^{-1})\left(t-T(\epsilon)+\frac{F(Mx(T(\epsilon)))}{\Lambda(\epsilon)}  \right).
\end{align*}
Since $T(\epsilon)>x_1(\epsilon)$, we have $\Lambda(\epsilon)T(\epsilon)>\Lambda(\epsilon)x_1(\epsilon)=F(Mx(T(\epsilon)))$. Therefore
$c:=T(\epsilon)-F(Mx(T(\epsilon)))/\Lambda(\epsilon)>0$. Hence for $t\geq T(\epsilon)$, we have
\begin{equation*}
(1-\Lambda(\epsilon))(1-\epsilon)\phi(x_U(t))
\geq
h(\Lambda(\epsilon))(1-\epsilon/2)(1-\epsilon)  (\phi\circ F^{-1})(t-c).
\end{equation*}
We may also write $T(\epsilon)=c+x_1(\epsilon)$. Therefore, for $t\geq T(\epsilon)$, we have
$t>t-c\geq T(\epsilon)-c=x_1(\epsilon)$. By \eqref{eq.x1Fdel}, for $t\geq T(\epsilon)$, we have
$t>t-c>x_1(\epsilon)>F(\delta)$. Therefore for $t\geq T(\epsilon)$, $F^{-1}(t)<F^{-1}(t-c)<\delta$, so
as $\phi$ is increasing on $(0,\delta)$ we have $(\phi\circ F^{-1})(t)<(\phi\circ F^{-1})(t-c)$. Hence for $t>T(\epsilon)$ we have
\begin{align*}
(1-\Lambda(\epsilon))(1-\epsilon)\phi(x_U(t))&>
h(\Lambda(\epsilon))(1-\epsilon/2)(1-\epsilon)  (\phi\circ F^{-1})(t)\\
&=L\frac{1+\epsilon}{(1-\epsilon)^2} (1-\epsilon/2)(1-\epsilon)  (\phi\circ F^{-1})(t)\\
&>L(1+\epsilon) (\phi\circ F^{-1})(t).
\end{align*}
Since $T(\epsilon)>T_1(\epsilon)$, we have $g(t)<L(1+\epsilon)(\phi\circ F^{-1})(t)$ for $t>T(\epsilon)$, so
\[
(1-\Lambda(\epsilon))(1-\epsilon)\phi(x_U(t))>g(t), \quad t\geq T(\epsilon).
\]
This implies
\[
x_U'(t)>-f(x_U(t))+g(t), \quad t\geq T(\epsilon); \quad x_U(T(\epsilon))>x(T(\epsilon)).
\]
A comparison argument now confirms that $x(t)<x_U(t)$ for all $t\geq T(\epsilon)$. By the definition of $x_U$ and the fact that $F$ is decreasing,
we have
\[
F(x(t))>\Lambda(\epsilon)(t-T(\epsilon))+F(Mx(T(\epsilon))), \quad t\geq T(\epsilon).
\]
Therefore
\[
\liminf_{t\to\infty} \frac{F(x(t))}{t}\geq \Lambda(\epsilon).
\]
Since $\Lambda(\epsilon)\to \Lambda_\ast$ as $\epsilon\to 0^+$, we get
\begin{equation} \label{eq.liminfFxt}
\liminf_{t\to\infty} \frac{F(x(t))}{t}\geq \Lambda_\ast.
\end{equation}

We now construct a lower solution, and obtain a companion limit superior bound to \eqref{eq.liminfFxt}.
From \eqref{eq.liminfFxt}, for every $\epsilon\in (0,1)$, there is a $T_1(\epsilon)>0$ such that $F(x(t))\geq \Lambda_\ast(1-\epsilon)$
for all $t\geq T_1(\epsilon)$. Hence $x(t)\leq F^{-1}(\Lambda_\ast(1-\epsilon))$ for $t\geq T_1(\epsilon)$. Since $F^{-1}\in \text{RV}_\infty(-1/(\beta-1))$,
we have
\begin{align*}
\limsup_{t\to\infty} \frac{x(t)}{F^{-1}(t)}&\leq \limsup_{t\to\infty} \frac{F^{-1}(\Lambda_\ast(1-\epsilon))}{F^{-1}(t)}\\
&=\lim_{t\to\infty} \frac{F^{-1}(\Lambda_\ast(1-\epsilon))}{F^{-1}(t)}
=(\Lambda_\ast(1-\epsilon))^{-1/(\beta-1)}.
\end{align*}
Therefore
\[
\limsup_{t\to\infty} \frac{x(t)}{F^{-1}(t)}\leq \Lambda_\ast^{-1/(\beta-1)}.
\]
Hence for any $\lambda>0$ we have
\begin{align}
\limsup_{t\to\infty} \frac{x(t)}{F^{-1}(\lambda t)}
&=\limsup_{t\to\infty} \frac{x(t)}{F^{-1}(t)}\cdot \frac{F^{-1}(t)}{F^{-1}(\lambda t)}\nonumber\\
\label{eq.xFinvlamt}
&\leq \Lambda_\ast^{-1/(\beta-1)} \cdot \lambda^{1/(\beta-1)}
=\left(\frac{\lambda}{\Lambda^\ast}\right)^{1/(\beta-1)}.
\end{align}
Define $\lambda(\epsilon)\in (0,1)$ by
\[
h(\lambda(\epsilon))=L\frac{1-\epsilon}{(1+\epsilon)(1+\sqrt{\epsilon})}<L=h(\Lambda_\ast).
\]
Since $h$ is decreasing, we have $\lambda(\epsilon)>\Lambda_\ast$. Thus
\[
\left(\frac{\lambda}{\Lambda^\ast}\right)^{1/(\beta-1)}>1.
\]
Also we have that $\lambda(\epsilon)\downarrow \Lambda_\ast$ as $\epsilon\to 0^+$.
By \eqref{eq.xFinvlamt}, there exists $T_1'(\epsilon)>0$ such that
\[
\frac{x(t)}{F^{-1}(\lambda(\epsilon) t)}\leq 2 \left(\frac{\lambda(\epsilon)}{\Lambda^\ast}\right)^{1/(\beta-1)}=:\frac{M}{2}<M, \quad t>T_1'(\epsilon).
\]
Clearly $M=4\left(\frac{\lambda(\epsilon)}{\Lambda^\ast}\right)^{1/(\beta-1)}>4$.

For every $\epsilon\in (0,1)$, there is $T_2'(\epsilon)>0$ such that
\[
g(t)>L(1-\epsilon)(\phi\circ F^{-1})(t), \quad t\geq T_2'(\epsilon).
\]
Define $T_3'(\delta,\epsilon)=F(\delta)/\lambda(\epsilon)$. Also we have that there is $x_3(\epsilon)>0$ such that
\[
\frac{f(x)}{\phi(x)}<1+\epsilon, \quad x<x_3(\epsilon).
\]
Since $x(t)\to 0$ as $t\to\infty$, it follows that there is a $T_4'(\epsilon)>0$ such that $x(t)<x_3(\epsilon)$ for all $t\geq T_4'(\epsilon)$.
Since $\phi\circ F^{-1}\in \text{RV}_\infty(-\beta/(\beta-1))$, for every $\eta\in (\epsilon,1)$, there exists $x_2'(\eta)>0$ such that
\[
\frac{(\phi\circ F^{-1})(\lambda(\epsilon)x)}{(\phi\circ F^{-1})(x)}\leq \lambda(\epsilon)^{-\beta/(\beta-1)}(1+\eta), \quad x\geq x_2'(\eta).
\]
Fix $\eta=\sqrt{\epsilon}$. Then
\[
\frac{(\phi\circ F^{-1})(\lambda(\epsilon)x)}{(\phi\circ F^{-1})(x)}\leq \lambda(\epsilon)^{-\beta/(\beta-1)}\left(1+\sqrt{\epsilon}\right), \quad
x\geq x_2'(\sqrt{\epsilon})=:x_4'(\epsilon).
\]
Finally, there is $T_5'(\delta)>0$ such that $x(t)<\delta$ for all $t\geq T_5'(\delta)$.
Define $T'(\epsilon)=1+\max(T_1'(\epsilon),x_4'(\epsilon),T_2'(\epsilon),T_3'(\delta,\epsilon),T_5'(\delta),T_4'(\epsilon))$ and
\[
x_L(t)=F^{-1}\left(\lambda(\epsilon)(t-T'(\epsilon))+F\left(\tfrac{x(T'(\epsilon))}{M}\right)\right), \quad t\geq T'(\epsilon).
\]
Then $x_L(T'(\epsilon))=x(T'(\epsilon))/M<x(T'(\epsilon))$ because $M>1$. Also as $T'(\epsilon)>T_5'(\delta)$, we have that
$x_L(t)<x(T'(\epsilon))<\delta$ for all $t\geq T'(\epsilon)$. For $t\geq T'(\epsilon)$, we have
$F'(x_L(t))x_L'(t)=\lambda(\epsilon)$. Therefore we have
\[
x_L'(t)=-f(x_L(t))+(1-\lambda(\epsilon))f(x_L(t)), \quad t\geq T'(\epsilon).
\]
Since $f(x)<(1+\epsilon)\phi(x)$ for $x<x_3(\epsilon)$, and $x_L(t)<x(T'(\epsilon))<x_3(\epsilon)$ for all $t>T'(\epsilon)$
we have $f(x_L(t))<(1+\epsilon)\phi(x_L(t))$ so
\[
x_L'(t)<-f(x_L(t))+(1-\lambda(\epsilon))(1+\epsilon)\phi(x_L(t)), \quad t\geq T'(\epsilon).
\]
Define $c'=F(\frac{1}{M}x(T'(\epsilon)))-\lambda(\epsilon)T'(\epsilon)$. Since $T'(\epsilon)>T_1(\epsilon)$, we have 
$x(T'(\epsilon))<MF^{-1}(\lambda(\epsilon)T'(\epsilon))$, so as $F$ is decreasing, it follows that  
$F(\frac{1}{M}x(T'(\epsilon)))>\lambda(\epsilon)T'(\epsilon)$ or $c'>0$.
Hence 
\begin{align*}
(1-\lambda(\epsilon))(1+\epsilon)\phi(x_L(t))
&=
(1-\lambda(\epsilon))(1+\epsilon)(\phi\circ F^{-1})(\lambda(\epsilon)t +c'). 
\end{align*}
Since $\lambda(\epsilon)T_3'(\delta,\epsilon)=F(\delta)$, for $t\geq T'(\epsilon)>T_3'(\delta,\epsilon)$ we have $\lambda(\epsilon)t>F(\delta)$. 
Thus $F^{-1}(\lambda(\epsilon)t)<\delta$. Since $c'>0$ and $F^{-1}$ is decreasing, it follows that  
$F^{-1}(\lambda(\epsilon)t+c')<F^{-1}(\lambda(\epsilon)t)<\delta$ for $t\geq T'(\epsilon)$, and as $\phi$ is increasing on $(0,\delta)$ we deduce that  
\[
(\phi\circ F^{-1})(\lambda(\epsilon)t+c')<(\phi\circ F^{-1})(\lambda(\epsilon)t), \quad t\geq T'(\epsilon).
\] 
Hence
\[
(1-\lambda(\epsilon))(1+\epsilon)\phi(x_L(t))<(1-\lambda(\epsilon))(1+\epsilon)(\phi\circ F^{-1})(\lambda(\epsilon)t), \quad t\geq T'(\epsilon),
\]
and so
\[
x_L'(t)<-f(x_L(t))+(1-\lambda(\epsilon))(1+\epsilon)(\phi\circ F^{-1})(\lambda(\epsilon)t), \quad t\geq T'(\epsilon).
\]
Since $t\geq T'(\epsilon)>T_4'(\epsilon)$, we have that 
\[
\frac{(\phi\circ F^{-1})(\lambda(\epsilon)t)}{(\phi\circ F^{-1})(t)}\leq \lambda(\epsilon)^{-\beta/(\beta-1)}\left(1+\sqrt{\epsilon}\right),
\]
so by the definition of $h$, we obtain
\[
x_L'(t)<-f(x_L(t))+h(\lambda(\epsilon))
(1+\epsilon)\left(1+\sqrt{\epsilon}\right)(\phi\circ F^{-1})(t), 
\quad t\geq T'(\epsilon).
\]
Therefore as 
$t\geq T'(\epsilon)>T_2'(\epsilon)$, we have 
$g(t)>L(1-\epsilon)(\phi\circ F^{-1})(t)$ for $t\geq T'(\epsilon)$, so 
for $t\geq T'(\epsilon)$ by using the definition of $\lambda(\epsilon)$, we get  
\begin{align*}
x_L'(t)
&<-f(x_L(t))+\frac{1}{L(1-\epsilon)} h(\lambda(\epsilon))(1+\epsilon)\left(1+\sqrt{\epsilon}\right) \cdot g(t)\\
&=-f(x_L(t))+ g(t).
\end{align*}
Since $x_L(T'(\epsilon))<x(T'(\epsilon))$, it follows that $x_L(t)<x(t)$ as $t\geq T'(\epsilon)$. Since $F$ is decreasing, using the definition of $x_L$ 
we arrive at 
\[
F(x(t))<\lambda(\epsilon)(t-T'(\epsilon))+F(\frac{1}{M}x(T'(\epsilon))), \quad t\geq T'(\epsilon).
\]
Therefore 
\[
\limsup_{t\to\infty} \frac{F(x(t))}{t}\leq \lambda(\epsilon).
\]
Letting $\epsilon\to 0^+$, and recalling that $\lambda(\epsilon)\to \Lambda_\ast$ as $\epsilon\to 0^+$ we get
\[
\limsup_{t\to\infty} \frac{F(x(t))}{t}\leq \Lambda_\ast.
\]
Combining this with \eqref{eq.liminfFxt} yields \eqref{eq.limFxtLambda} as required. 

\subsection{Proof of Theorem~\ref{theorem.Fxtto0}}
We consider first the proof when $\beta>1$, and then sketch the proof when 
$f\circ F^{-1}\in \text{RV}_\infty(-1)$ and $F^{-1}\in \text{RV}_\infty(0)$. Let $\epsilon\in (0,1/2)$. 
Since $f(x)/\phi(x)\to1$ as $x\to 0^+$, we have that there exists $x_1(\epsilon)>0$ such that $f(x)<(1+\epsilon)\phi(x)$ for all $x\leq x_1(\epsilon)$. 
Since $f\in \text{RV}_0(\beta)$, it follows that $\phi\in \text{RV}_0(\beta)$ and therefore that $h:=\phi\circ F^{-1}\in \text{RV}_\infty(-\beta/(\beta-1))$.
By \eqref{eq.gfFinfty}, we have that there exists $T_1(\epsilon)>0$ such 
that $h(t)<\epsilon^{1+\beta/(\beta-1)} g(t)$ for $t\geq T_1(\epsilon)$. Also, as $h\in \text{RV}_\infty(-\beta/(\beta-1))$, we have that $h(\epsilon t)/h(t)\to \epsilon^{-\beta/(\beta-1)}$ as $t\to\infty$. 
Hence there exists $T_2(\epsilon)>0$ such that $h(\epsilon t)<2\epsilon^{-\beta/(\beta-1)}h(t)$ for $t\geq T_2(\epsilon)$. Define $T(\epsilon)=1+\max(T_1(\epsilon),T_2(\epsilon))$ and 
\begin{equation} \label{def.xslowM}
M=\max\left(2,\frac{2x(T)}{x_1(\epsilon)},\frac{x(T)}{F^{-1}(\epsilon T)}\right).
\end{equation}
Also define
\begin{equation} \label{def.xslowxL}
x_L(t)=F^{-1}\biggl(\epsilon(t-T)+F(x(T)/M)\biggr), \quad t\geq T.
\end{equation}
Since $M>1$, we have $x_L(T)=x(T)/M<x(T)$. Also, the definition of $M$ implies that $x(T)/M\leq x_1(\epsilon)/2$. Since $x_L$ is decreasing on $[T,\infty)$ 
it follows that $x_L(t)\leq x(T)/M<x_1(\epsilon)$ for $t\geq T$. Hence for $t\geq T$,
\[
f(x_L(t))<(1+\epsilon)\phi(x_L(t))=(1+\epsilon)h(\epsilon(t-T)+F^\ast),
\]
where $F^\ast:=F(x(T))/M$. Since $M\geq x(T)/F^{-1}(\epsilon T)$ and $F$ is decreasing, we have $\epsilon T\leq F^\ast$. Therefore for $t\geq T$ 
we have $\epsilon(t-T)+F^\ast\geq \epsilon t$. Since $h$ is decreasing, we have $h(\epsilon t)\geq h(\epsilon(t-T)+F^\ast)$. Therefore 
\[
f(x_L(t))<(1+\epsilon)h(\epsilon t), \quad t\geq T.
\]
We note by definition that $x_L$ is in $C^1(T,\infty)$ and $x_L'(t)=-\epsilon f(x_L(t))$ for $t\geq T$. Then for $t\geq T$, we get
\[
x_L'(t)=-f(x_L(t))+(1-\epsilon)f(x_L(t))<-f(x_L(t))+(1-\epsilon^2) h(\epsilon t)<-f(x_L(t))+h(\epsilon t).
\]
Let $t\geq T$. Since $T>T_2$, we have that $h(\epsilon t)<2\epsilon^{-\beta/(\beta-1)}h(t)$ and as $T>T_1$ we have $h(t)<\epsilon^{1+\beta/(\beta-1)} g(t)$. 
Therefore as $\epsilon<1/2$ we have
\begin{align*}
x_L'(t)&<-f(x_L(t))+h(\epsilon t)< -f(x_L(t)) + 2\epsilon^{-\beta/(\beta-1)}h(t)< -f(x_L(t))+  2\epsilon g(t) \\ 
&<-f(x_L(t))+g(t)
\end{align*}
for $t\geq T$. Therefore we have that $x_L(t)<x(t)$ for $t\geq T$. Since $x(t)\to 0$ as $t\to\infty$, and there exists $\delta_1>0$ such that $F(x)>0$ for 
all $x<\delta_1$, it follows that there is $T'>0$ such that $x(t)<\delta_1$ for $t\geq T'$ and therefore $F(x(t))>0$ for $t\geq T'$. Hence for $t\geq \max(T',T)$ 
we have $0<F(x(t))<F(x_L(t))=\epsilon(t-T)+F^\ast$. Therefore, we get 
\[
0\leq \liminf_{t\to\infty} \frac{F(x(t))}{t}\leq \limsup_{t\to\infty} \frac{F(x(t))}{t}\leq \epsilon.
\] 
Letting $\epsilon\to 0^+$ finally gives \eqref{eq.xslow}, as required.

Suppose now that $f\circ F^{-1}\in \text{RV}_\infty(-1)$ and $F^{-1}\in \text{RV}_\infty(0)$. Since $f(x)/\phi(x)\to 1$ as $x\to 0^+$ and $F^{-1}(t)\to 0$ 
as $t\to\infty$, we have
\[
\lim_{t\to\infty} \frac{\phi(F^{-1}(t))}{f(F^{-1}(t))}=1.
\]
Hence $h=\phi\circ F^{-1}\in \text{RV}_\infty(-1)$. Let $\epsilon\in (0,1/2)$. By \eqref{eq.gfFinfty}, we have that there exists $T_1(\epsilon)>0$ such 
that $h(t)<\epsilon^2 g(t)$ for $t\geq T_1(\epsilon)$. Also, as $h\in \text{RV}_\infty(-1)$, we have that $h(\epsilon t)/h(t)\to \epsilon^{-1}$ as $t\to\infty$. 
Hence there exists $T_2(\epsilon)>0$ such that $h(\epsilon t)<2\epsilon^{-1}h(t)$ for $t\geq T_2(\epsilon)$. Define $T(\epsilon)=1+\max(T_1(\epsilon),T_2(\epsilon))$. Now define $M$ and $x_L$ as in \eqref{def.xslowM} and \eqref{def.xslowxL}. Proceeding in a manner identical to that used in the case when $\beta>1$, we can show that 
once again that 
\[
x_L'(t)<-f(x_L(t))+g(t), \quad t\geq T; \quad x_L(T)<x(T).
\]
Therefore, we have that $x_L(t)<x(t)$ for $t\geq T$, and proceeding as in the case when $\beta>1$, it can once more be shown that \eqref{eq.xslow} holds.

\subsection{Proof of Theorem~\ref{theorem.fxtgt1}}
We prove part (i). Since $f\in \text{RV}_0(\beta)$ there is an increasing $\varphi\in C^1(0,\infty)$ such that 
\begin{equation} \label{eq.fphiRV}
\lim_{x\to 0^+} \frac{f(x)}{\varphi(x)}=1, \quad \lim_{x\to 0^+} \frac{x\varphi'(x)}{\varphi(x)}=\beta.
\end{equation}
Since $g \in \text{RV}_\infty(-\theta)$ and $\theta>0$, there exists a decreasing $\gamma\in C^1(0,\infty)$ such that 
\begin{equation} \label{eq.ggammaRV}
\lim_{t\to\infty} \frac{g(t)}{\gamma(t)}=1, \quad \lim_{t\to\infty} \frac{t\gamma'(t)}{\gamma(t)}=-\theta.
\end{equation}
\textbf{STEP 1: An estimate deriving from \eqref{eq.gfFinfty}:}
We will show that \eqref{eq.gfFinfty} implies 
\begin{equation} \label{eq.phiinvgammatgamma0}
\lim_{t\to\infty} \frac{\varphi^{-1}(\gamma(t))}{t\gamma(t)}=0.
\end{equation}
Since $g(t)/\gamma(t)\to 1$ as $t\to\infty$, we have from \eqref{eq.gfFinfty} implies $f(F^{-1}(t))/\gamma(t)\to 0$ as $t\to\infty$. Thus for every $\epsilon\in(0,1)$, there exists $T(\epsilon)>0$ such that $f(F^{-1}(t))<\epsilon\gamma(t)$ for all $t\geq T(\epsilon)$ or $\epsilon^{-1}f(F^{-1}(t))<\gamma(t)$
for all $t\geq T(\epsilon)$. Since $\gamma$ is decreasing, $\gamma^{-1}$ is also, and therefore $\gamma^{-1}(\epsilon^{-1}f(F^{-1}(t)))>t$ for $t\geq T(\epsilon)$. Thus
\[
\frac{t}{\gamma^{-1}(f\circ F^{-1})(t)}<\frac{\gamma^{-1}(\frac{1}{\epsilon}(f\circ F^{-1})(t))}{\gamma^{-1}((f\circ F^{-1})(t))}, \quad t\geq T(\epsilon).
\] 
Since $f\in\text{RV}_0(\beta)$, it follows that $f\circ F^{-1}\in \text{RV}_\infty(-\beta/(\beta-1))$. Thus $(f\circ F^{-1})(t)\to 0$ as $t\to\infty$. 
Also, as $\gamma\in \text{RV}_\infty(-\theta)$ and $\gamma$ is decreasing, we have that $\gamma^{-1}\in \text{RV}_0(-1/\theta)$. Therefore 
\[
\lim_{x\to 0^+} \frac{\gamma^{-1}(\epsilon^{-1}x)}{\gamma^{-1}(x)}=\epsilon^{1/\theta}.
\]
Therefore 
\[
\limsup_{t\to\infty} \frac{t}{\gamma^{-1}(f\circ F^{-1})(t)}\leq \epsilon^{1/\theta}.
\]
Since $\epsilon>0$ is arbitrary, we have 
\[
\lim_{t\to\infty} \frac{t}{\gamma^{-1}(f\circ F^{-1})(t)}=0,
\]
or 
\[
\lim_{t\to\infty} \frac{\gamma^{-1}(f\circ F^{-1})(t)}{t}=+\infty.
\]
Since $F^{-1}(t)\to 0$ as $t\to\infty$, we have $f(F^{-1}(t))/\varphi(F^{-1}(t))\to 1$ as $t\to\infty$, and therefore, because $\gamma^{-1}\in \text{RV}_0(-1/\theta)$, we have 
\[
\lim_{t\to\infty} \frac{\gamma^{-1}(f(F^{-1}(t)))}{\gamma^{-1}(\varphi(F^{-1}(t)))}=1.
\]
Hence
\[
\lim_{t\to\infty} \frac{\gamma^{-1}((\varphi\circ F^{-1}))(t)}{t}=+\infty.
\]
Since $(\varphi\circ F^{-1})(t)\to 0$ as $t\to\infty$ and $\varphi\circ F^{-1}$ is invertible with inverse $F\circ \varphi^{-1}$ we get
\[
\lim_{x\to 0^+} \frac{\gamma^{-1}(x)}{(F\circ \varphi^{-1})(x)}=+\infty.
\]
Since $f\in \text{RV}_0(\beta)$ and $\beta>1$ we have that 
\[
\lim_{x\to 0^+} \frac{F(x)}{x/f(x)}=\frac{1}{\beta-1},
\]
and hence
\[
\lim_{x\to 0^+} \frac{F(x)}{x/\varphi(x)}=\frac{1}{\beta-1}.
\]
Since $\varphi^{-1}(y)\to 0$ as $y\to 0^+$ we have 
\[
\lim_{y\to 0^+} \frac{(F\circ \varphi^{-1})(y)}{\varphi^{-1}(y)/y}=\frac{1}{\beta-1}.
\]
Hence 
\[
\lim_{x\to 0^+} \frac{\gamma^{-1}(x)}{\varphi^{-1}(x)/x}
=\lim_{x\to 0^+} 
\frac{(F\circ \varphi^{-1})(x)}{\varphi^{-1}(x)/x}\cdot 
\frac{\gamma^{-1}(x)}{(F\circ \varphi^{-1})(x)}=+\infty.
\]
Since $\gamma(t)\to 0$ as $t\to\infty$, we arrive at 
\[
\lim_{t\to\infty} \frac{t}{\varphi^{-1}(\gamma(t))/\gamma(t)} =+\infty,
\]
which implies \eqref{eq.phiinvgammatgamma0}, as required. 

\textbf{STEP 2: Lower bound.} We determine a lower bound on the solution. Let $\epsilon\in (0,1)$. Then, there exists $T_1(\epsilon)>0$ such that 
\[
\frac{g(t)}{\gamma(t)}>1-\epsilon, \quad t\geq T_1(\epsilon).
\]
Also, for every $\epsilon\in (0,1/2)$ there exists $x_1(\epsilon)>0$ such that 
\[
\frac{f(x)}{\varphi(x)}\leq 1+\epsilon, \quad x\leq x_1(\epsilon).
\]
Define $T_2(\epsilon)=\sup\{t>0\,:\, x(t)\leq x_1(\epsilon)/2\}$, 
$T(\epsilon)=1+\max(T_1(\epsilon),T_2(\epsilon))$, and 
\[
K_\epsilon:=\gamma^{-1}\left(\varphi(2^{-1}x(T(\epsilon)))\right)>0.
\]
Finally, define 
\[
x_L(t)=\varphi^{-1}\left( \frac{1-\epsilon}{1+\epsilon}\gamma(t+K_\epsilon)  \right), \quad t\geq T(\epsilon).
\]
Then as $\varphi^{-1}$ is increasing, and $\gamma$ is decreasing we have 
\begin{align*}
x_L(T(\epsilon))
&=\varphi^{-1}\left( \frac{1-\epsilon}{1+\epsilon}\gamma(T(\epsilon)+K_\epsilon)  \right)
<\varphi^{-1}\left( \gamma(T(\epsilon)+K_\epsilon)  \right)\\
&<\varphi^{-1}\left( \gamma(K_\epsilon)  \right)
=\varphi^{-1}\left( \varphi(x(T(\epsilon)/2)) \right)=x(T(\epsilon))/2\\
&<x(T(\epsilon)).
\end{align*}
For $t\geq T(\epsilon)$ we have that $x_L(t)\leq x_L(T(\epsilon))<x(T(\epsilon))$ and because $T(\epsilon)>T_1(\epsilon)$
we get $x(T(\epsilon))\leq x_1(\epsilon)/2<x_1(\epsilon)$. Therefore for $t\geq T(\epsilon)$, we have $x_L(t)<x_1(\epsilon)$. 
Thus
\[
\frac{f(x_L(t))}{\varphi(x_L(t))}\leq 1+ \epsilon, \quad t\geq T(\epsilon). 
\]
Hence by the definition of $x_L$, and the fact that $\gamma$ is decreasing, we have for $t\geq T(\epsilon)$
\[
f(x_L(t))\leq (1+\epsilon)\varphi(x_L(t)) = (1-\epsilon)\gamma(t+K_\epsilon)<(1-\epsilon)\gamma(t)<g(t).
\]
Since $\varphi\in C^1(0,\infty)$, $\varphi'(x)>0$ for all $x>0$, $\gamma\in C^1(0,\infty)$ and $\gamma'(t)<0$ for all $t>0$, we see 
that $x_L'(t)$ is well--defined for all $t>T(\epsilon)$ and is given by 
\[
x_L'(t)=\frac{1}{\varphi'(x_L(t))}\cdot \frac{1-\epsilon}{1+\epsilon} \gamma'(t+K_\epsilon)<0.
\]
Therefore we have 
\[
x_L'(t)<-f(x_L(t))+g(t), \quad t\geq T(\epsilon),
\]
and since $x_L(T(\epsilon))<x(T(\epsilon))$, we have that $x_L(t)<x(t)$ for all $t\geq T(\epsilon)$. Since $\varphi$ is increasing, we have 
$\varphi(x_L(t))<\varphi(x(t))$ for all $t\geq T(\epsilon)$. Hence
\[
\varphi(x(t))> \frac{1-\epsilon}{1+\epsilon}\gamma(t+K_\epsilon), \quad t\geq T(\epsilon).
\]
Since $\gamma\in \text{RV}_\infty(-\theta)$, it follows that $\gamma(t+K_\epsilon)/\gamma(t)\to 1$ as $t\to\infty$. Therefore 
\[
\liminf_{t\to\infty} \frac{\varphi(x(t))}{\gamma(t)}\geq \frac{1-\epsilon}{1+\epsilon}.
\]
Letting $\epsilon\to 0^+$ and recalling that $f(x)/\varphi(x)\to 1$ as $x\to 0^+$ and $g(t)/\gamma(t)\to 1$ as $t\to\infty$, 
we arrive at 
\begin{equation} \label{eq.fxtgtliminf}
\liminf_{t\to\infty} \frac{f(x(t))}{g(t)}\geq 1.
\end{equation}
\textbf{STEP 3: Upper bound.} We need a limit superior to companion \eqref{eq.fxtgtliminf}. 
For every $\epsilon\in (0,1/2)$, there exists $x_0(\epsilon)>0$ such that 
\[
1-\epsilon<\frac{f(x)}{\varphi(x)}<1+\epsilon, \quad x<x_0(\epsilon).
\]
Let $M_\epsilon:=1+\epsilon>1$. Since $x(t)\to 0$ as $t\to\infty$, and $\varphi(x)\to 0$ as $x\to 0^+$, 
it follows that there is $T_0(\epsilon)>0$ such that 
\[
M_\epsilon\varphi(x(t))\leq \frac{1-\epsilon}{2}\cdot \frac{\varphi(x_0(\epsilon))}{1+2\epsilon}, \quad t\geq T_0(\epsilon). 
\] 
Since $g(t)/\gamma(t)\to 1$ as $t\to\infty$, there exists $T_1(\epsilon)>0$ such that 
\[
\frac{g(t)}{\gamma(t)}<1+\epsilon, \quad t\geq T_1(\epsilon).
\]
By \eqref{eq.ggammaRV}, there also exists $T_2(\epsilon)>0$ such that 
\[
\frac{-t\gamma'(t)}{\gamma(t)}<\theta+\epsilon, \quad t\geq T_2(\epsilon).
\]
Furthermore, by \textbf{STEP 2}, there exists a $T_3(\epsilon)>0$ such that 
\[
\frac{\varphi(x(t))}{\gamma(t)}>\frac{1}{1+\epsilon}=\frac{1}{M_\epsilon}, \quad t\geq T_3(\epsilon).
\]
By \eqref{eq.fphiRV}, we have that $x\varphi'(x)/\varphi(x)\to \beta$ as $x\to 0^+$. It therefore follows that there exists $x_1>0$ such that 
\[
\frac{\varphi(x)}{x\varphi'(x)}<\frac{2}{\beta}, \quad x<x_1.
\]
Since $\varphi(x(t))\to 0$ as $t\to\infty$ there exists $T_4(\epsilon)>0$ such that
\[
M_\epsilon\varphi(x(t))\leq \frac{1}{2}\gamma(T_2(\epsilon)), \quad t\geq T_4(\epsilon),
\]
and $T_5(\epsilon)>0$ such that 
\[
M_\epsilon \varphi(x(t))\leq \frac{1}{2}\frac{1-\epsilon}{1+2\epsilon}\varphi(x_1), \quad t\geq T_5(\epsilon).
\]
Since $\varphi^{-1}\in \text{RV}_0(1/\beta)$, there exists $x_2(\epsilon)>0$ such that 
\[
\frac{\varphi^{-1}\left(\frac{1+2\epsilon}{1-\epsilon}x  \right)}{\varphi^{-1}(x)}<2\left(\frac{1+2\epsilon}{1-\epsilon}\right)^{1/\beta}, \quad x<x_2(\epsilon).
\]
Since $\varphi(x(t))\to 0$ as $t\to\infty$, there exists $T_6(\epsilon)>0$ such that 
\[
 M_\epsilon \varphi(x(t))\leq x_2(\epsilon)/2, \quad t\geq T_6(\epsilon). 
\]
Since \eqref{eq.phiinvgammatgamma0} holds, it follows that there is $T_7(\epsilon)>0$ such that
\[
\frac{\varphi^{-1}(\gamma(t))}{t\gamma(t)}<\frac{1}{2}\frac{\epsilon}{4^{1/\beta}(4\theta+2)/\beta}, \quad t\geq T_7(\epsilon),
\]
and because $\varphi(x(t))\to 0$ as $t\to\infty$ we have that there exists $T_8(\epsilon)>0$ such that 
\[
M_\epsilon \varphi(x(t))\leq \gamma(T_7(\epsilon))/2, \quad t\geq T_8(\epsilon). 
\] 
Now define 
\[
T(\epsilon)=1+\max(T_0(\epsilon), T_1(\epsilon), T_2(\epsilon), T_3(\epsilon), T_4(\epsilon), T_5(\epsilon), T_6(\epsilon), T_7(\epsilon), T_8(\epsilon)),
\]
and with $\varphi_1(\epsilon):=\varphi(x(T(\epsilon)))$, we further define  
\[
x_U(t)=\varphi^{-1}\biggl(\frac{1+2\epsilon}{1-\epsilon} \gamma\bigl(t-T(\epsilon)+\gamma^{-1}(M_\epsilon \varphi_1(\epsilon))\bigr) \biggr), \quad t\geq T(\epsilon).
\]
Since $\varphi^{-1}$ is increasing and $M_\epsilon>1$, it follows that 
\[
x_U(T(\epsilon))=\varphi^{-1}\left( \frac{1+2\epsilon}{1-\epsilon} \varphi(x(T(\epsilon))) \right)>\varphi^{-1}(\varphi(x(T(\epsilon))))=x(T(\epsilon)).
\]
Since $T(\epsilon)>T_0(\epsilon)$, have that 
\[
M_\epsilon\varphi(x(T(\epsilon)))\leq \frac{1-\epsilon}{2}\cdot \frac{\varphi(x_0(\epsilon))}{1+2\epsilon}
<\frac{1-\epsilon}{1+2\epsilon} \varphi(x_0(\epsilon)).
\]
Therefore
\[
\varphi^{-1}\left(\frac{1+2\epsilon}{1-\epsilon}M_\epsilon\varphi(x(T(\epsilon)))\right)< x_0(\epsilon),
\]
so $x_U(T(\epsilon))<x_0(\epsilon)$. Hence $x_U(t)<x_0(\epsilon)$ for all $t\geq T(\epsilon)$. 

We start by estimating $-f(x_U(t))+g(t)$ for all $t\geq T(\epsilon)$. Since $x_U(t)<x_0(\epsilon)$ for all $t\geq T(\epsilon)$, we have 
\[
-(1-\epsilon)\varphi(x_U(t))>-f(x_U(t))>-(1+\epsilon)\varphi(x_U(t)), \quad t\geq T(\epsilon).
\]
Define $\varphi_1(\epsilon):=\varphi(x(T(\epsilon)))$ and
\[
c:=T(\epsilon)-\gamma^{-1}(M_\epsilon\varphi_1(\epsilon)).
\] 
Since $t\geq T(\epsilon)>T_1(\epsilon)$, by the definition of $x_U$, we have 
\[
-f(x_U(t))+g(t)<-(1-\epsilon)\varphi(x_U(t))+(1+\epsilon)\gamma(t)
=-(1+2\epsilon)\gamma(t-c)+(1+\epsilon)\gamma(t).
\]
Since $t\geq T(\epsilon)>T_3(\epsilon)$, we have $\varphi_1(\epsilon)=\varphi(x(T(\epsilon)))>\gamma(T(\epsilon))/M_\epsilon$.
Hence $\gamma^{-1}(M_\epsilon\varphi_1(\epsilon))<T(\epsilon)$. Thus $c=T(\epsilon)-\gamma^{-1}(M_\epsilon\varphi_1(\epsilon))>0$. Moreover
$T(\epsilon)=c+\gamma^{-1}(M_\epsilon\varphi_1(\epsilon))$. Hence for $t\geq T(\epsilon)$, we have $t-c>0$ and so
$\gamma(t-c)>\gamma(t)$ for all $t\geq T(\epsilon)$. Therefore for all $t\geq T(\epsilon)$, we have 
\begin{multline} \label{eq.rhsupper1}
-f(x_U(t))+g(t)<-(1+2\epsilon)\gamma(t-c)+(1+\epsilon)\gamma(t)\\
<-(1+2\epsilon)\gamma(t-c)+(1+\epsilon)\gamma(t-c)=-\epsilon \gamma(t-c).
\end{multline}

We now seek to estimate $x_U'(t)$ for all $t\geq T(\epsilon)$. Since $(\varphi^{-1})'(x)=1/\varphi'(\varphi^{-1}(x))$, for $t\geq T(\epsilon)$ we have 
\[
-x_U'(t)=-\frac{1}{\varphi'(x_U(t))}\cdot \frac{1+2\epsilon}{1-\epsilon}\gamma'(t-c),
\]
or 
\[
-x_U'(t)=\frac{1+2\epsilon}{1-\epsilon} \cdot \frac{\varphi(x_U(t))}{x_U(t)\varphi'(x_U(t))}\cdot\frac{x_U(t)}{t-c} \cdot\frac{1}{\varphi(x_U(t))}\cdot   \frac{-(t-c)\gamma'(t-c)}{\gamma(t-c)}\cdot\gamma(t-c).
\]
Since $T(\epsilon)>T_5(\epsilon)$, we have that 
\[
M_\epsilon\varphi(x(T(\epsilon)))\leq \frac{1}{2}\frac{1-\epsilon}{1+2\epsilon}\varphi(x_1).
\]
Also 
\[
\varphi(x_U(T(\epsilon)))=\frac{1+2\epsilon}{1-\epsilon}M_\epsilon\varphi(x(T(\epsilon)))
\leq  \frac{1}{2}\varphi(x_1)<\varphi(x_1),
\]
so $x_U(T(\epsilon))<x_1$. Hence $x_U(t)<x_1$ for all $t\geq T(\epsilon)$. Hence
\[
\frac{\varphi(x_U(t))}{x_U(t)\varphi'(x_U(t))}<\frac{2}{\beta}, \quad t\geq T(\epsilon).  
\]
Since $\varphi(x_U(t))=(1+2\epsilon)/(1-\epsilon) \gamma(t-c)$ for $t\geq T(\epsilon)$, we arrive at
\[
-x_U'(t)< \frac{2}{\beta} \cdot\frac{\varphi^{-1}\left((1+2\epsilon)/(1-\epsilon) \gamma(t-c)\right)}{(t-c)\gamma(t-c)}\cdot   \frac{-(t-c)\gamma'(t-c)}{\gamma(t-c)}\cdot\gamma(t-c).
\]
Now for $t\geq T(\epsilon)$, we have $t-c\geq T(\epsilon)-c=\gamma^{-1}(M_\epsilon\varphi_1(\epsilon))$. Since $T(\epsilon)>T_4(\epsilon)\geq T_2(\epsilon)$
we have $\varphi(M_\epsilon \varphi(x(T(\epsilon))))\leq \gamma(T_2(\epsilon))/2$. Hence 
$M_\epsilon \varphi_1(\epsilon)<\gamma(T_2(\epsilon))$. Hence
$\gamma^{-1}(M_\epsilon\varphi_1(\epsilon))>T_2(\epsilon)$. Thus for $t\geq T(\epsilon)$, we have 
$t-c\geq \gamma^{-1}(M_\epsilon\varphi_1(\epsilon))>T_2(\epsilon)$. Hence 
\[
-\frac{(t-c)\gamma'(t-c)}{\gamma(t-c)}<\theta+\epsilon, \quad t\geq T(\epsilon). 
\] 
Thus
\[
-x_U'(t)< \frac{2(\theta+\epsilon)}{\beta} \cdot\frac{\varphi^{-1}\left((1+2\epsilon)/(1-\epsilon) \gamma(t-c)\right)}{(t-c)\gamma(t-c)}\cdot\gamma(t-c).
\]
Next, we note for $t\geq T(\epsilon)$ that $t-c\geq \gamma^{-1}(M_\epsilon\varphi_1(\epsilon))$. Thus $\gamma(t-c)\leq M_\epsilon\varphi_1(\epsilon)$. 
Since $T(\epsilon)>T_6(\epsilon)$, we have that $M_\epsilon\varphi(x(T(\epsilon)))\leq x_2(\epsilon)/2<x_2(\epsilon)$, or 
$M_\epsilon\varphi_1(\epsilon))<x_2(\epsilon)$. Thus for all $t\geq T(\epsilon)$ we have  $\gamma(t-c)< x_2(\epsilon)$. Hence 
\[
\frac{\varphi^{-1}\left(\frac{1+2\epsilon}{1-\epsilon}\gamma(t-c)\right)}{\varphi^{-1}(\gamma(t-c))}<2\left(\frac{1+2\epsilon}{1-\epsilon}\right)^{1/\beta}.
\]
Hence for $t\geq T(\epsilon)$, we have 
\[
-x_U'(t)< \frac{2(\theta+\epsilon)}{\beta} 
\cdot\frac{2\left(\frac{1+2\epsilon}{1-\epsilon}\right)^{1/\beta} \varphi^{-1}(\gamma(t-c))}{(t-c)\gamma(t-c)}
\cdot\gamma(t-c),
\]
and so as $\epsilon<1/2$, we have 
\begin{align*}
-x_U'(t)&< \frac{4(\theta+\epsilon)}{\beta}  \left(\frac{1+2\epsilon}{1-\epsilon}\right)^{1/\beta}
\cdot\frac{\varphi^{-1}(\gamma(t-c))}{(t-c)\gamma(t-c)} \cdot \gamma(t-c)\\
&<\frac{4\theta+2}{\beta} 4 ^{1/\beta} \cdot\frac{\varphi^{-1}(\gamma(t-c))}{(t-c)\gamma(t-c)} \cdot \gamma(t-c).
\end{align*}
Since $t\geq T(\epsilon)$, we have $t-c\geq T(\epsilon)-c=\gamma^{-1}(M_\epsilon\varphi_1(\epsilon))$. Since $t\geq T(\epsilon)>T_8(\epsilon)$, 
we have that $M_\epsilon \varphi(x(T(\epsilon)))\leq \gamma(T_7(\epsilon))/2<\gamma(T_7(\epsilon))$. Hence $\gamma^{-1}(M_\epsilon \varphi_1(\epsilon))>T_7(\epsilon)$, and so $t-c>T_7(\epsilon)$ for all $t\geq T(\epsilon)$. This yields
\[
\frac{\varphi^{-1}(\gamma(t-c))}{(t-c)\gamma(t-c)}<\frac{1}{2}\cdot \frac{\epsilon}{\frac{4\theta+2}{\beta} 4^{1/\beta}}, \quad t\geq T(\epsilon).
\] 
Hence for $t\geq T(\epsilon)$ we arrive at 
\[
-x_U'(t)<\frac{4\theta+2}{\beta} 4 ^{1/\beta} \cdot\frac{\varphi^{-1}(\gamma(t-c))}{(t-c)\gamma(t-c)} \cdot \gamma(t-c)
\frac{1}{2} \epsilon \gamma(t-c).
\]
Using this inequality and \eqref{eq.rhsupper1}, it follows for $t\geq T(\epsilon)$ that 
\[
x_U'(t)> -\frac{1}{2}\epsilon\gamma(t-c)>-\epsilon\gamma(t-c)>-f(x_U(t))+g(t).
\]
Since $x_U(T(\epsilon))>x(T(\epsilon))$, we have that $x_U(t)>x(t)$ for all $t\geq T(\epsilon)$. Therefore by the definition of $x_U$, we get
\[
\varphi(x(t)) < \varphi(x_U(t))=\frac{1+2\epsilon}{1-\epsilon}\gamma(t-c), \quad t\geq T(\epsilon).
\]
Since $\gamma\in \text{RV}_\infty(-\theta)$, it follows that $\gamma(t-c)/\gamma(t)\to 1$ as $t\to\infty$. Therefore 
\[
\limsup_{t\to\infty} \frac{\varphi(x(t))}{\gamma(t)}\leq \frac{1+2\epsilon}{1-\epsilon}.
\]
Letting $\epsilon\to 0^+$, and recalling that $x(t)\to 0$ as $t\to\infty$ and that $f(x)/\varphi(x)\to 1$ as $x\to 0^+$, we get 
\[
\limsup_{t\to\infty} \frac{f(x(t))}{\gamma(t)}\leq 1.
\]
Since $\gamma(t)/g(t)\to 1$ as $t\to\infty$, we have
\[
\limsup_{t\to\infty} \frac{f(x(t))}{g(t)}\leq 1.
\]
Combining this limit with \eqref{eq.fxtgtliminf} yields \eqref{eq.fxtgtlim}, as required.

We now turn to the proof of part (ii). It is assumed that $g$ is asymptotic to a decreasing function; this function must be positive. Call it $\gamma_1$. The fact that $g$ is slowly varying implies that there exists a positive function $\gamma_2\in C^1(0,\infty)$ such that $g(t)/\gamma_2(t)\to 1$ as $t\to\infty$ and $t\gamma_2'(t)/\gamma_2(t)\to 0$ as $t\to\infty$. We note that the estimate \eqref{eq.phiinvgammatgamma0} derived in \textbf{STEP 1} holds, with $\gamma_1$ in the role of $\gamma$ even though the proof above does not apply in the case $\theta=0$. In fact, in this slowly varying case, the proof is much easier. Since $\gamma_1\in \text{RV}_\infty(0)$ and $\varphi\in \text{RV}_0(\beta)$, we have that $\varphi^{-1}\circ \gamma_1\in \text{RV}_\infty(0)$. On the other hand, $t\mapsto t\gamma_1(t)$ is a function in $\text{RV}_\infty(1)$. Therefore we have that $t\mapsto \varphi^{-1}(\gamma_1(t))/(t\gamma_1(t))$ is in $\text{RV}_\infty(-1)$. Thus $\varphi^{-1}(\gamma_1(t))/(t\gamma_1(t))\to 0$ as $t\to\infty$. 

We start by establishing a lower bound on the solution of \eqref{eq.odepert}, just as in \textbf{STEP 2} above. Since $g(t)/\gamma_1(t)\to 1$ and $g(t)/\gamma_2(t)\to 1$ as $t\to\infty$, we have that $\gamma_2(t)/\gamma_1(t)\to 1$ as $t\to\infty$, and so $\varphi^{-1}(\gamma_2(t))/\varphi^{-1}(\gamma_1(t))\to 1$ as $t\to\infty$. Hence 
\[
\lim_{t\to\infty} \frac{\varphi^{-1}(\gamma_2(t))}{t\gamma_2(t)}=0.
\]
We construct the following estimates. Let $\epsilon\in (0,1/4)$. Then there exists $T_1(\epsilon)>0$ such that 
\[
\frac{g(t)}{\gamma_1(t)}>\frac{1-\epsilon}{1+\epsilon/2}, \quad 1-\epsilon<\frac{\gamma_2(t)}{\gamma_1(t)}<1+\epsilon, \quad t\geq T_1(\epsilon).
\]
Since $x\varphi'(x)/\varphi(x)\to \beta$ as $x\to 0^+$ and $f(x)/\varphi(x)\to 1$ as $x\to 0^+$, for every $\epsilon\in (0,1/4)$ there exists $x_1(\epsilon)>0$
such that 
\[
\frac{f(x)}{\varphi(x)}<1+\epsilon, \quad x<x_1(\epsilon)
\]
and there is an $x_2>0$ such that 
\[
\frac{\varphi(x)}{x\varphi'(x)}<\frac{2}{\beta}, \quad x< x_2.
\]
Since $\varphi^{-1}(\gamma_2(t))/(t\gamma_2(t))\to 0$ as $t\to\infty$, we have that there exists $T_2(\epsilon)>0$ such that 
\[
\frac{\varphi^{-1}(\gamma_2(t))}{t\gamma_2(t)}<\frac{\sqrt{\epsilon}}{20/\beta}, \quad t\geq T_2(\epsilon).
\]
Since $\gamma_1(t)\to 0$ as $t\to\infty$, there exists $T_3(\epsilon)>0$ such that 
\[
\gamma_1(t)<\varphi(x_2), \quad \gamma(t)<\varphi(x_1(\epsilon)), \quad t\geq T_3(\epsilon).
\]
Since $\gamma_2$ obeys $t\gamma_2'(t)/\gamma_2(t)\to 0$ as $t\to\infty$, it follows that there exists $T_4(\epsilon)> 0$ such that 
\[
\left|\frac{t\gamma_2'(t)}{\gamma_2(t)}\right|<\sqrt{\epsilon}, \quad t\geq T_4(\epsilon).
\]
Define $T(\epsilon)=1+\max(T_1,T_2,T_3,T_4)$. Finally, define 
\[
K=\gamma_1^{-1}\left(\varphi(x(T)/2)\right)>0,
\]
(so $\gamma_1(K)=\varphi(x(T)/2)$) and 
\[
x_L(t)=\varphi^{-1}\left(\frac{1-\epsilon}{(1+\epsilon)^3}\gamma_2(t+K)\right), \quad t\geq T.
\]
For $t\geq T$, we have $\gamma_2(t+K)<(1+\epsilon)\gamma_1(t+K)<(1+\epsilon)\gamma_1(K)$. Thus as $\varphi^{-1}$ is increasing, we have
\[
x_L(T)<\varphi^{-1}\left(\frac{1-\epsilon}{(1+\epsilon)^2}\gamma_1(K)\right)<\varphi^{-1}\left(\gamma_1(K)\right)=x(T)/2<x(T).
\]
For $t\geq T>T_3$, we also have 
\[
x_L(t)<\varphi^{-1}\left(\frac{1-\epsilon}{(1+\epsilon)^2}\gamma_1(t+K)\right)<\varphi^{-1}(\gamma_1(T))<x_1(\epsilon),
\]
and also
\[
x_L(t)<\varphi^{-1}(\gamma_1(T))<x_2.
\]
Since $K>0$, for $t\geq T$ by the definition of $x_L$ and the monotonicity of $\gamma_1$ we have 
\[
f(x_L(t))<(1+\epsilon)\varphi(x_L(t))=\frac{1-\epsilon}{(1+\epsilon)^2}\gamma_2(t+K)<\frac{1-\epsilon}{1+\epsilon}\gamma_1(t+K)<\frac{1-\epsilon}{1+\epsilon}\gamma_1(t).
\] 
On the other hand, for $t\geq T$ we have
\[
g(t)>\frac{1-\epsilon}{1+\epsilon/2}\gamma_1(t).
\]
Hence
\[
f(x_L(t))-g(t)< (1-\epsilon)\left\{\frac{1}{1+\epsilon}- \frac{1}{1+\epsilon/2}\right\}\gamma_1(t)=-\epsilon \gamma_1(t) 
\frac{1-\epsilon}{2(1+\epsilon/2)(1+\epsilon)}.
\]
Since $\epsilon\in (0,1/4)$, we have that 
\[
\frac{1-\epsilon}{2(1+\epsilon/2)(1+\epsilon)}>\frac{4}{15}.
\]
Therefore
\[
f(x_L(t))-g(t)< -\epsilon\cdot \frac{4}{15} \gamma_1(t), \quad t\geq T.
\]

On the other hand, since $x_L\in C^1(T,\infty)$ and has derivative given by 
\[
x_L'(t)\varphi(x_L(t))=\frac{1-\epsilon}{(1+\epsilon)^3}\gamma_2'(t+K)
\]
using the definition of $x_L$ we have for $t\geq T$
\begin{multline*}
x_L'(t)=\frac{\varphi(x_L(t))}{x_L(t)\varphi'(x_L(t))}\cdot \frac{\varphi^{-1}\left(\frac{1-\epsilon}{(1+\epsilon)^3}\gamma_2(t+K)\right)}{\varphi^{-1}(\gamma_2(t+K))}\\
\times \frac{\varphi^{-1}(\gamma_2(t+K))}{(t+K)\gamma_2(t+K)}\cdot \frac{(t+K)\gamma_2'(t+K)}{\gamma_2(t+K)}\cdot \gamma_2(t+K).
\end{multline*}
Since $x_L(t)<x_1(\epsilon)$ for $t\geq T$, the first factor is bounded above by $2/\beta$. Since $\varphi^{-1}$ is increasing, the second factor is bounded above by unity. Since $t+K>t\geq T$, the absolute value of the fourth factor is bounded by $\sqrt{\epsilon}$. Hence
\[
|x_L'(t)|\leq \frac{2}{\beta}\cdot \frac{\varphi^{-1}(\gamma_2(t+K))}{(t+K)\gamma_2(t+K)}\cdot \sqrt{\epsilon} \frac{\gamma_2(t+K)}{\gamma_1(t+K)}\cdot \gamma_1(t+K).
\]
and as $\gamma_2(t+K)<(1+\epsilon)\gamma_1(t+K)<5/4\cdot \gamma_1(t+K)$ and $\gamma_1$ is decreasing, we have
\[
|x_L'(t)|\leq \frac{\varphi^{-1}(\gamma_2(t+K))}{(t+K)\gamma_2(t+K)}\cdot \sqrt{\epsilon} \frac{5}{\beta} \cdot \gamma_1(t).
\]
Since $t+K>t\geq T$ we have 
\[
\frac{\varphi^{-1}(\gamma_2(t+K))}{(t+K)\gamma_2(t+K)} < \frac{\sqrt{\epsilon}}{20/\beta},
\]
so 
\[
|x_L'(t)|<  \frac{\epsilon}{4} \gamma_1(t).
\]
Hence for $t\geq T$ we have
\[
-x_L'(t)> -\frac{\epsilon}{4} \gamma_1(t)> -\epsilon\cdot \frac{4}{15} \gamma_1(t) > f(x_L(t))-g(t).
\]
Hence 
\[
x_L'(t)<-f(x_L(t))+g(t), \quad t\geq T; \quad x(T)<x_L(T).
\]
Thus $x_L(t)<x(t)$ for $t\geq T$. Hence $\varphi(x(t))>\varphi(x_L(t))=(1-\epsilon)(1+\epsilon)^{-3}\gamma_2(t+K)$. Therefore as $\gamma_2(t+K)/\gamma_2(t)\to 1$
as $t\to\infty$, we have 
\[
\liminf_{t\to\infty} \frac{\varphi(x(t))}{\gamma_2(t)}\geq (1-\epsilon)(1+\epsilon)^{-3}.
\]
Letting $\epsilon\to 0^+$ and noting that $f(x)/\varphi(x)\to 1$ as $x\to 0^+$ and $g(t)/\gamma_2(t)\to 1$ as $t\to\infty$, we have 
\[
\liminf_{t\to\infty} \frac{f(x(t))}{g(t)}\geq 1.
\]

To obtain the upper bound, set $m=4^{1+2/\beta}/\beta$, let $\epsilon\in (0,1/2)$ and define $M_\epsilon=1+\epsilon$. There exist $x_0$, $x_1$ and $x_2$ such that 
\begin{gather*}
1-\epsilon<\frac{f(x)}{\varphi(x)}<1+\epsilon, \quad x<x_0(\epsilon),\\
\frac{\varphi(x)}{x\varphi'(x)}<\frac{2}{\beta}, \quad x<x_1, \\
\frac{\varphi^{-1}\left(\frac{(1+2\epsilon)^2}{(1-\epsilon)^2}x\right)}{\varphi^{-1}(x)}<2\left(\frac{1+2\epsilon}{1-\epsilon}\right)^{2/\beta}, \quad x<x_2(\epsilon). 
\end{gather*}
For every $\epsilon\in (0,1/2)$, there exists a $T_1(\epsilon)$ such that 
\[
\frac{g(t)}{\gamma_1(t)}<(1+\epsilon)^2, \quad 1-\epsilon<\frac{\gamma_2(t)}{\gamma_1(t)}<1+\epsilon, \quad t\geq T_1(\epsilon).
\]
Similarly, for every $\epsilon\in (0,1/2)$, there exists a $T_2(\epsilon)$ such that 
\[
\left|\frac{t\gamma_2'(t)}{\gamma_2(t)}\right|<\sqrt{\epsilon}, \quad t\geq T_2(\epsilon).
\]
From the lower bound, there is $T_3(\epsilon)>0$ such that 
\[
\frac{\varphi(x(t))}{\gamma_1(t)}>\frac{1}{1+\epsilon}=\frac{1}{M}, \quad t\geq T_3(\epsilon). 
\]
Since $\varphi(x(t))\to 0$ as $t\to\infty$, we may define $T_4(\epsilon)>0$ such that 
\[
M\varphi(x(t))\leq \frac{1}{2}\gamma_1(T_1(\epsilon)), \quad M\varphi(x(t))\leq \frac{1}{2}\gamma_1(T_2(\epsilon)), \quad t\geq T_4(\epsilon),
\]
and similarly there is a $T_5(\epsilon)>0$ such that 
\[
\varphi(x(t))<\frac{1}{2}\left(\frac{1-\epsilon}{(1+\epsilon)(1+2\epsilon)} \right)^2\varphi(x_0\wedge x_1), \quad t\geq T_5(\epsilon),
\]
where, as is conventional, $a\wedge b$ denotes the minimum of the real numbers $a$ and $b$.
Furthermore, there exists $T_6(\epsilon)>0$ such that 
\[
\frac{\varphi^{-1}(\gamma_2(t))}{t\gamma_2(t)}<\frac{\sqrt{\epsilon}}{m}, \quad t\geq T_6(\epsilon). 
\]
We define $T_7(\epsilon)>0$ such that 
\[
M\varphi(x(t))\leq \frac{1}{2}\gamma_1(T_6(\epsilon)), \quad t\geq T_7(\epsilon).
\]
Since $\gamma_2(t)\to 0$ as $t\to\infty$, there is $T_8(\epsilon)>0$ such that 
\[
\gamma_2(t)\leq \frac{1}{2}x_2(\epsilon), \quad t\geq T_8(\epsilon),
\]
and finally there is a $T_9(\epsilon)>0$ such that 
\[
M\varphi(x(t))\leq \frac{1}{2}\gamma_1(T_8(\epsilon)), \quad t\geq T_9(\epsilon).
\]
Define $T=1+\max_{j=1,\ldots,9} T_j(\epsilon)$, $K:=\gamma_1^{-1}(M\varphi(x(T)))$ (suppressing here the $\epsilon$--dependence) and define 
\[
x_U(t)=\varphi^{-1}\left(\frac{(1+2\epsilon)^2}{(1-\epsilon)^2} \gamma_2(t-T+K)\right), \quad t\geq T.
\] 
Since $T>T_3(\epsilon)$, we have $\gamma_1(K)=M\varphi(x(T))>\gamma_1(T)$, and because $\gamma_1$ is decreasing, we have that $K<T$. Define $c:=T-K>0$. 

Since $T>T_4(\epsilon)$, $\gamma_1(K)=M\varphi(x(T))\leq \gamma_1(T_1(\epsilon))/2<\gamma_1(T_1(\epsilon))$, so $K>T_1(\epsilon)$. Thus $\gamma_2(K)>(1-\epsilon)\gamma_1(K)$. Hence as $\varphi^{-1}$
is increasing and $M>1$ we have
\[
x_U(T)=\varphi^{-1}\left(\frac{(1+2\epsilon)^2}{(1-\epsilon)^2} \gamma_2(K)\right)>\varphi^{-1}\left(\frac{(1+2\epsilon)^2}{1-\epsilon} \gamma_1(K)\right)
>\varphi^{-1}(\gamma_1(K))>x(T).
\]
For $t\geq T$, we have $t-T+K\geq K>T_1(\epsilon)$, so $\gamma_2(t-c)<(1+\epsilon)\gamma_1(t-T+K)<(1+\epsilon)\gamma_1(K)$. Thus for $t\geq T$ we have 
\[
x_U(t)<\varphi^{-1}\left(\frac{(1+2\epsilon)^2}{(1-\epsilon)^2}(1+\epsilon)\gamma_1(K) \right)
=\varphi^{-1}\left(\frac{(1+2\epsilon)^2(1+\epsilon)^2}{(1-\epsilon)^2}\varphi(x(T))\right).
\]
Since $T>T_5(\epsilon)$, the argument of $\varphi^{-1}$ is less than $\varphi(x_0(\epsilon))$ and also less than $\varphi(x_1)$. Therefore 
\[
x_U(t)<x_0(\epsilon), \quad x_U(t)<x_1, \quad t\geq T.
\] 
Hence for $t\geq T$, we have 
\[
f(x_U(t))>(1-\epsilon)\varphi(x_U(t))=\frac{(1+2\epsilon)^2}{1-\epsilon}\gamma_2(t-c).
\]
Since $t\geq T$, $t-c\geq T-c=K>T_1(\epsilon)$, we have $\gamma_2(t-c)>(1-\epsilon)\gamma_1(t-c)$. Hence for $t\geq T$, we have 
\[
f(x_U(t))>(1+2\epsilon)^2\gamma_1(t-c)>(1+2\epsilon)^2 \gamma_1(t)
\] 
as $c>0$ and $\gamma_1$ is decreasing. Since $t\geq T>T_1(\epsilon)$, we have $g(t)<(1+\epsilon)^2\gamma_1(t)$, so
\[
-f(x_U(t))+g(t)<\left\{-(1+2\epsilon)^2 +(1+\epsilon)^2\right\}\gamma_1(t)=-2\epsilon(1+3\epsilon/2)\gamma_1(t)<-2\epsilon\gamma_1(t).
\]
Next, we can use the definition of $x_U$ to obtain the identity
\begin{multline*}
x_U'(t)=\left(\frac{1+2\epsilon}{1-\epsilon}\right)^2 \frac{\varphi(x_U(t))}{x_U(t)\varphi'(x_U(t))}\frac{1}{\varphi(x_U(t))}\frac{(t-c)\gamma_2'(t-c)}{\gamma_2(t-c)}\frac{x_U(t)}{t-c}\gamma_2(t-c)
\end{multline*} 
for $t>T$. For $t\geq T$, $t-T\geq K=\gamma_1^{-1}(M\varphi(x(T)))$. Since $T>T_4(\epsilon)$, $M\varphi(x(T))<\gamma_1(T_2(\epsilon))$, so $K>T_2(\epsilon)$. Thus $t\geq T$ implies $t-c>T_2(\epsilon)$. 
Hence 
\[
\left|\frac{(t-c)\gamma_2'(t-c)}{\gamma_2(t-c)}\right|<\sqrt{\epsilon}, \quad t\geq T,
\]
and because $x_U(t)<x_1$ for $t\geq T$, we have 
\[
\frac{\varphi(x_U(t))}{x_U(t)\varphi'(x_U(t))}<\frac{2}{\beta}, \quad t\geq T.
\]
Using these estimates and the definition of $x_U$, we obtain 
\[
|x_U'(t)|<  \frac{2}{\beta} \sqrt{\epsilon}  \frac{1}{\gamma_2(t-c)} \frac{x_U(t)}{t-c}\gamma_2(t-c), \quad t\geq T,
\]
and again using the definition of $x_U$ we have 
\[
|x_U'(t)|<  \frac{2}{\beta} \sqrt{\epsilon}  \frac{\varphi^{-1}(\gamma_2(t-c))}{(t-c)\gamma_2(t-c)} \frac{\varphi^{-1}\left(\left(\frac{1+2\epsilon}{1-\epsilon}\right)^2\gamma_2(t-c)\right)}{\varphi^{-1}(\gamma_2(t-c))}\gamma_2(t-c), \quad t\geq T.
\]
For $t\geq T$, $t-c\geq K$. Since $t>T_7(\epsilon)$, $M\varphi(x(T))<\gamma_1(T_6(\epsilon))$, so $K>T_6(\epsilon)$. Hence $t\geq T$ implies 
$t-c>T_6(\epsilon)$,  we obtain  
\[
\frac{\varphi^{-1}(\gamma_2(t-c))}{(t-c)\gamma_2(t-c)}<\frac{\sqrt{\epsilon}}{m}, \quad t\geq T.
\]
For $t\geq T$, $t-c\geq K$. Since $t>T_9(\epsilon)$, $M\varphi(x(T))<\gamma_1(T_8(\epsilon))$, so $K>T_8(\epsilon)$. 
Hence $t\geq T$ implies $t-c>T_8(\epsilon)$, so we have 
$\gamma_2(t-c)\leq x_2(\epsilon)/2<x_2(\epsilon)$ and therefore
\[
\frac{\varphi^{-1}\left(\left(\frac{1+2\epsilon}{1-\epsilon}\right)^2\gamma_2(t-c)\right)}{\varphi^{-1}(\gamma_2(t-c))}<2 \left(\frac{1+2\epsilon}{1-\epsilon} \right)^{2/\beta}, \quad t\geq T.
\]
Using these estimates, we get
\[
|x_U'(t)|<  \frac{4}{\beta m} \epsilon  \left(\frac{1+2\epsilon}{1-\epsilon} \right)^{2/\beta} \gamma_2(t-c), \quad t\geq T.
\]
Since $\epsilon<1/2$, $c>0$, $\gamma_1$ is decreasing, by the definition of $m$, we get
\[
|x_U'(t)|<  \frac{4}{\beta m} \epsilon  4^{2/\beta} \frac{\gamma_2(t-c)}{\gamma_1(t-c)} \cdot \gamma_1(t)
=\epsilon   \frac{\gamma_2(t-c)}{\gamma_1(t-c)} \cdot \gamma_1(t), \quad t\geq T.
\]
For $t\geq T$ we have $t-c\geq K$. Since $T>T_4(\epsilon)$, $M\varphi(x(T))<\gamma_1(T_1(\epsilon))$, so $K>T_1(\epsilon)$. Hence $t-c>T_1(\epsilon)$, which implies $\gamma_2(t-c)/\gamma_1(t-c)<1+\epsilon<3/2$ for $t\geq T$. Hence
\[
|x_U'(t)|<  \frac{3}{2}\epsilon \gamma_1(t), \quad t\geq T,
\]
so $x_U'(t)>-\frac{3}{2}\gamma_1(t)$ for $t\geq T$. Recall that $-f(x_U(t))+g(t)<-2\epsilon \gamma_1(t)$ for $t\geq T$. Hence 
\[
x_U'(t)>-\frac{3}{2}\gamma_1(t)>-2\epsilon \gamma_1(t)>-f(x_U(t))+g(t), \quad t\geq T,
\]
and since $x_U(T)>x(T)$, it follows that $x_U(t)>x(t)$ for $t\geq T$. Hence 
\[
\varphi(x(t))<\varphi(x_U(t))=\left(\frac{1+2\epsilon}{1-\epsilon}\right)^2 \gamma_2(t-c), \quad t\geq T. 
\]
Since $\gamma_2\in \text{RV}_\infty(0)$, we have 
\[
\limsup_{t\to\infty} \frac{\varphi(x(t))}{\gamma_2(t)}\leq \left(\frac{1+2\epsilon}{1-\epsilon}\right)^2,
\]
and letting $\epsilon\to 0^+$ yields
\[
\limsup_{t\to\infty} \frac{\varphi(x(t))}{\gamma_2(t)}\leq 1.
\]
Since $x(t)\to 0$ as $t\to\infty$, we get
\[
\limsup_{t\to\infty} \frac{f(x(t))}{g(t)}\leq 1.
\]
Combining this with the lower estimate gives $f(x(t))/g(t)\to 1$ as $t\to\infty$ as required. 

\begin{flushleft}
\textbf{Acknowledgements}
\end{flushleft}
The authors are pleased to thank the organisers of the conference for the opportunity to present their work at the meeting. 
We are also extremely grateful for the very careful checking of the calculations and keen scrutiny of the paper by the referee. 
We finally wish to thank Robert Szczelina (Jagiellonian University) for valuable technical advice concerning this submission.

\end{document}